\documentclass[fullpage]{article}
\RequirePackage[OT1]{fontenc}
\RequirePackage{color}
\RequirePackage{amsthm,amsmath}
\RequirePackage[numbers]{natbib}
\RequirePackage[colorlinks,citecolor=blue,urlcolor=blue]{hyperref}
\RequirePackage{hypernat}
\def\sDiag{{\mathop{\hbox{\scriptsize\rm Diag}}}}
\def\bA{{\mathbf{A}}}
\def\bG{{\mathbf{G}}}
\def\inter{{\mathop{\hbox{\rm int}\,}}}
\def\argmin{\mathop{\hbox{\rm argmin\,}}}

\def\Opt{{\mathop{\hbox{\rm Opt}}}}
\def\bR{{\mathbf{R}}}
\def\errvi{{\mathop{\hbox{\rm Err}_{\hbox{\rm vi}}}}}
\def\errNash{{\mathop{\hbox{\rm Err}_{\hbox{\scriptsize\rm N}}}}}
\def\bS{{\mathbf{S}}}
\def\Prob{{\mathop{\hbox{\rm Prob}}}}
\def\Tr{\mathop{\hbox{\rm Tr}}}
\def\Diag{\mathop{\hbox{\rm Diag}}}
\def\rint{\hbox{\rm rint}}

\def\cX{{\cal X}}
\def\cY{{\cal Y}}
\def\cA{{\cal A}}

\def\cS{{\cal S}}
\def\cF{{\cal F}}
\def\cH{{\cal H}}

\def\cB{{\cal B}}
\def\cC{{\cal C}}

\def\cE{{\cal E}}

\def\bS{{\mathbf{S}}}

\def\bE{{\mathbf{E}}}
\def\Prob{{\hbox{\rm Prob}}}
\newtheorem{remark}{Remark}
\newtheorem{algorithm}{Algorithm}
\newtheorem{theorem}{Theorem}
\newtheorem{corollary}{Corollary}
\newtheorem{lemma}{Lemma}
\newtheorem{proposition}{Proposition}

\def\argmin{\mathop{\rm argmin}}

\newcommand{\g} {\gamma}
\newcommand{\wh}[1]{\widehat{#1}}

\newcommand{\e}{\epsilon}

\newcommand{\cD}{{\cal D}}

\newcommand{\cG}{{\cal G}}

\newcommand{\bP}{{\bf P}}

\newcommand{\be}{\begin{eqnarray}}
\newcommand{\ee}[1]{\label{eq:#1}\end{eqnarray}}
\newcommand{\nn}{\nonumber \\}
\newcommand{\ese}{\end{eqnarray*}}
\newcommand{\bse}{\begin{eqnarray*}}
\newcommand{\rf}[1]{~(\ref{eq:#1})}

\newtheorem{theo}{Theorem}
\newcommand{\bthm}{\begin{theo}}
\newcommand{\ethm}[1]{\label{the:#1}\end{theo}\par}

\newtheorem{col}{Corollary}
\newcommand{\bcol}{\begin{col}}
\newcommand{\ecol}[1]{\label{the:#1}\end{col}\par}

\newtheorem{defi}{Definition}
\newcommand{\bdf}{\begin{defi}}
\newcommand{\edf}[1]{\label{df:#1}\end{defi}\par}

\newtheorem{lem}{Lemma}
\newcommand{\blem}{\begin{lem}}
\newcommand{\elem}[1]{\label{le:#1}\end{lem}\par}

\newtheorem{pro}{Proposition}
\newcommand{\bpro}{\begin{pro}}
\newcommand{\epro}[1]{\label{pr:#1}\end{pro}\par}

\newtheorem{statement}{Problem}
\newcommand{\bstatement}{\begin{statement}}
\newcommand{\estatement}[1]{\label{stat:#1}\end{statement}\par}

\newcounter{assc}
\newcommand{\bass}[1]{\refstepcounter{assc}\label{ass:#1} \begin{maliste}{\bf \arabic{assc}.}}
\newcommand{\eass}{\end{maliste}}

\newcounter{algoc}
\newcommand{\balgo}[1]{\refstepcounter{algoc}\label{algo:#1} \begin{malistea}{\bf \arabic{algoc}.}}
\newcommand{\ealgo}{\end{malistea}}

\newcommand{\pr} {\noindent{\bf Proof\,:~}}
\newcommand{\epr}{\hfill\hbox{\hskip 4pt
                \vrule width 5pt height 6pt depth 1.5pt}\vspace{0.5cm}\par}

\title{Solving variational inequalities with Stochastic Mirror-Prox algorithm}
\author{Anatoli Juditsky\thanks{LJK,
Universit\'e J. Fourier, B.P. 53, 38041 Grenoble Cedex 9, France,
{\tt anatoli.juditsky@imag.fr}} \and Arkadi
Nemirovski\thanks{Georgia Institute
 of Technology, Atlanta, Georgia
30332, USA, {\tt nemirovs@isye.gatech.edu},
\newline Research of
this author was partly supported by the Office of Naval Research grant \# N000140811104 and the NSF grant DMS-0914785.}
\and Claire Tauvel\thanks{LJK,
Universit\'e J. Fourier, B.P. 53, 38041 Grenoble Cedex 9, France,
{\tt claire.tauvel@imag.fr}} }
\begin{document}
\thispagestyle{empty}
\maketitle
\begin{abstract}
In this paper we consider iterative methods for {\em stochastic
variational inequalities} (s.v.i.) with monotone operators. Our basic assumption is that the operator possesses both smooth and nonsmooth components. Further, only noisy observations of the problem data are available.
We develop a novel {Stochastic Mirror-Prox} (SMP)
algorithm for solving s.v.i. and show that with the convenient stepsize strategy it attains the optimal rates of convergence with respect to the problem parameters.
We apply the SMP algorithm to Stochastic composite minimization and describe particular applications to Stochastic Semidefinite Feasibility problem and deterministic Eigenvalue minimization.
\end{abstract}
\section{Introduction}
Variational inequalities with monotone operators form a convenient framework for unified treatment (including algorithmic design) of problems with ``convex structure,'' like convex minimization, convex-concave saddle point problems and convex Nash equilibrium problems. In this paper we utilize this framework to develop first order algorithms for {\sl stochastic} versions of the outlined problems, where the precise first order information is replaced with its unbiased stochastic estimates. This situation arises naturally in convex Stochastic Programming, where the precise first order information is unavailable (see examples in section 4). In some situations, e.g. those considered in \cite[Section 3.3]{Shapiroetal} and in Section 4.4,  where passing from available, but relatively computationally expensive precise first order information to its cheap stochastic estimates allows to accelerate the solution process, with the gain from randomization growing progressively with problem's sizes. \par
Our ``unifying framework'' is as follows.
Let $Z$ be a convex
compact set in Euclidean space $\cE$ with inner product
$\langle\cdot,\cdot\rangle$, $\|\cdot\|$ be a norm on $E$ (not
necessarily the one associated with the inner product), and
$F:Z\to \cE$ be a  monotone mapping: \be \forall(z,z'\in Z):\langle
F(z)-F(z'),z-z'\rangle\geq0\ee{F0} We are interested to
approximate a solution to the variational inequality (v.i.)
\begin{equation}\label{VI}
\hbox{find\ }z_*\in Z: \langle F(z),z_*-z\rangle \leq 0 \quad
\forall z\in Z
\end{equation}
associated with $Z,F$. Note that since $F$ is monotone on $Z$, the
condition in (\ref{VI}) is implied by  $\langle
F(z_*),z-z_*\rangle \geq0$ for all $z\in Z$, which is the standard
definition of a (strong) solution to the v.i. associated with
$Z,F$. The inverse -- a solution to v.i. as defined by (\ref{VI})
(a ``weak'' solution) is a strong solution as well -- also is true,
provided, e.g., that $F$ is continuous. An advantage of the
concept of weak solution is that such a solution  always exists under our
assumptions ($F$ is well defined and monotone on a convex compact
set $Z$).
\par
We quantify the inaccuracy of a candidate solution $z\in Z$ by
the error
\be
\errvi(z):=\max_{u\in Z}\langle F(u),z-u\rangle;
\ee{error}
note that this error  is always $\geq0$ and equals zero iff $z$ is a solution to (\ref{VI}).
\par
In what follows we impose on $F$, aside of the monotonicity, the
requirement
\be \forall(z,z'\in Z):\|F(z)-F(z')\|_*\leq
L\|z-z'\|+M \ee{F} with some known constants $L\geq0,M\geq0$. From now on,
\be\|\xi\|_*=\max\limits_{z:\|z\|\leq1}\langle \xi,z\rangle
\ee{normstar}
 is the norm
conjugate to $\|\cdot\|$.
\par  We are interested in the case
where (\ref{VI}) is solved by an iterative algorithm based on a
{\sl stochastic oracle representation} of  the operator
$F(\cdot)$. Specifically, when solving the problem, the algorithm
acquires information on $F$ via subsequent calls to a black box
 (``stochastic oracle'', SO). At the $i$th call, $i=0,1,...$, the
oracle gets as input a search point $z_i\in Z$ (this point is
generated by the algorithm on the basis of the information
accumulated so far) and returns the vector $\Xi(z_i,\zeta_i)$,
where $\{\zeta_i\in\bR^N\}_{i=1}^\infty$ is a sequence of i.i.d.
(and independent of the queries of the algorithm) random
variables. We suppose that the Borel function $\Xi(z,\zeta)$ is such that
\be
\forall z\in Z: \;\bE\left\{\Xi(z,\zeta_1)\right\}=F(z),\;\;
\bE\left\{\|\Xi(z,\zeta_i)-F(z)\|_*^2\right\}\leq \sigma^2.
\ee{oracle}
We call a monotone v.i. (\ref{eq:F0}), augmented by a stochastic oracle (SO), a {\sl stochastic} monotone v.i. (s.v.i.).
\par
To motivate our goal, let us start with known results \cite{NemYu} on the limits of performance
of iterative algorithms  for solving large-scale stochastic monotone v.i.'s.
To ``normalize'' the situation, assume that $Z$ is the unit Euclidean ball in $\cE=\bR^n$ and that $n$ is large. In this case,
the accuracy after $t$ steps of any algorithm for solving v.i.'s cannot be better than $O(1)\left[{L\over t}+{M+\sigma\over\sqrt{t}}\right]$.
In other words, for a properly chosen positive absolute constant $C$, for every number of steps $t$, all large enough values of $n$ and any algorithm $\cB$ for solving s.v.i.'s on the unit ball of $\bR^n$,  one can point out a monotone s.v.i. satisfying \rf{F}, \rf{oracle} and such that the expected error of the approximate solution $\tilde{z}_t$ generated by $\cB$ after $t$ steps , applied to such s.v.i.,
is  at least $c\left[{L\over t}+{M+\sigma\over\sqrt{t}}\right]$ for some $c>0$. To the best of our knowledge, no existing algorithm allows to achieve, uniformly in the dimension, this convergence rate. In fact, the ``best approximations'' available are given by Robust Stochastic Approximation (see \cite{Shapiroetal} and references therein) with the guaranteed rate of convergence $O(1){L+M+\sigma\over\sqrt{t}}$ and
extra-gradient-type algorithms for solving deterministic monotone v.i.'s with Lipschitz continuous operators (see \cite{MP,nes1,nes2,nes3}), which attain the accuracy $O(1){L\over t}$ in the case of  $M=\sigma=0$ or $O(1){M\over \sqrt{t}}$ when $L=\sigma=0$.
\par
The goal of this paper is to demonstrate that a specific {\sl Mirror-Prox} algorithm \cite{MP}
for solving monotone v.i.'s with Lipschitz continuous operators can be extended onto monotone s.v.i.'s to yield, uniformly in the dimension, the optimal rate of convergence $O(1)\left[{L\over t}+{M+\sigma\over\sqrt{t}}\right]$. We present the corresponding extension and investigate it in detail: we show how the algorithm can be ``tuned'' to the geometry of the s.v.i. in question and derive bounds for the probability of large deviations of the resulting error. We also present a number of applications where the specific structure of the rate of convergence indeed ``makes a difference.''
\par
The main body of the paper is organized as follows: in Section \ref{sect:Nash}, we
describe several special cases of monotone v.i.'s we are especially interested in (convex Nash equilibria, convex-concave saddle point problems, convex minimization). We single out these special cases since here one can define a useful ``functional'' counterpart $\errNash(\cdot)$ of the just defined error $\errvi(\cdot)$; both $\errNash$ and $\errvi$ will participate in our subsequent efficiency estimates. Our main development -- the {\sl Stochastic Mirror Prox} (SMP) algorithm -- is presented
in Section \ref{sect:SMP}. where we also provide some general results about its performance. Then in Section  \ref{SAa} we present SMP for Stochastic composite minimization and discuss its applications to Stochastic Semidefinite Feasibility problem and Eigenvalue minimization. All technical proofs are collected in the appendix.
\paragraph{Notations.} In the sequel, lowercase Latin letters denote vectors (and sometimes matrices). Script capital letters, like $\cE$, $\cY$, denote Euclidean spaces; the inner product in such a space, say, $\cE$, is denoted by $\langle\cdot,\cdot\rangle_\cE$ (or merely $\langle\cdot,\cdot\rangle$, when the corresponding space is clear from the context). Linear mappings from one Euclidean space to another, say, from $\cE$ to $\cF$, are denoted by boldface capitals like $\mathbf{A}$
(there are also some reserved boldface capitals, like $\bE$ for expectation, $\bR^k$ for the $k$-dimensional coordinate space, and $\bS^k$ for the space of $k\times k$ symmetric matrices).
$\mathbf{A}^*$ stands for the conjugate to mapping $\mathbf{A}$: if $\mathbf{A}: \cE\to\cF$, then $\mathbf{A}^*:\cF\to\cE$ is given by the identity $\langle f,\mathbf{A}e\rangle_\cF=\langle \mathbf{A}^*f,e\rangle_\cE$ for $f\in \cF,e\in\cE$. When both the origin and the destination space of a linear map, like $\mathbf{A}$, are the standard coordinate spaces, the map is identified with its matrix $A$, and $\mathbf{A}^*$ is identified with $A^T$. For a norm $\|\cdot\|$ on $\cE$, $\|\cdot\|_*$ stands for the conjugate norm, see \rf{normstar}. \par For Euclidean spaces $\cE_1,...,\cE_m$, $\cE=\cE_1\times...\times \cE_m$ denotes their Euclidean direct product, so that a vector from $\cE$ is a collection $u=[u_1;...;u_m]$ (``MATLAB notation'') of vectors $u_\ell\in\cE_\ell$, and $\langle u,v\rangle_\cE={\sum}_\ell\langle u_\ell,v_\ell\rangle_{\cE_\ell}$. Sometimes we allow ourselves to write $(u_1,...,u_m)$ instead of $[u_1;...;u_m]$.

\section{Preliminaries and Problem of interest} \label{sect:Nash}
\subsection{Nash v.i.'s and functional error}
In the sequel, we shall be especially interested in a special
case of v.i. (\ref{VI}) -- in a {\sl Nash} v.i. coming from a {\sl convex Nash Equilibrium} problem, and in the associated
{\sl functional} error measure. The Nash Equilibrium problem
can be described as follows: there are $m$ players, the $i$th of them choosing a point $z_i$ from a given set $Z_i$. The loss of the $i$th player is a given function $\phi_i(z)$ of the collection $z=(z_1,...,z_m)\in Z=Z_1\times...\times Z_m$ of players' choices.
With slight abuse of notation, we use for $\phi_i(z)$
also the notation $\phi_i(z_i,z^i)$, where $z^i$ is the collection of choices of all but the $i$th players.
Players are interested to minimize their losses, and Nash equilibrium $\widehat{z}$ is a point from $Z$ such that for every $i$ the function $\phi_i(z_i,\widehat{z}^i)$ attains its minimum in $z_i\in Z_i$ at $z_i=\widehat{z}_i$ (so that in the state $\widehat{z}$ no player has an incentive to change his choice, provided that the other players stick to their choices).\par
We call a Nash equilibrium problem {\sl convex}, if for every $i$, $Z_i$ is a compact convex set, $\phi_i(z_i,z^i)$ is a Lipschitz continuous function convex in $z_i$ and concave in $z^i$, and the function $\Phi(z)={\sum}_{i=1}^m\phi_i(z)$ is convex.
It is well known (see, e.g., \cite{NOR}) that setting
$$
F(z)=\left[F^1(z);\dots;F^m(z)\right],\,F^i(z)\in\partial_{z_i}\phi_i(z_i,z^i),\,i=1,...,m
$$
where $\partial_{z_i}\phi_i(z_i,z^i)$ is the subdifferential of
the convex function $\phi_i(\cdot,z^i)$ at a point $z_i$, we get a
monotone operator such that the solutions to the corresponding
v.i. (\ref{VI}) are exactly the Nash equilibria.
 Note that since $\phi_i$ are Lipschitz continuous, the associated operator
 $F$ can be chosen to be bounded. For this v.i. one can consider, along with
 the v.i.-accuracy measure $\errvi(z)$, the {\sl functional} error measure
$$
\errNash(z)={\sum}_{i=1}^m\left[\phi_i(z)-\min_{w_i\in
Z_i}\phi_i(w_i,z^i)\right]
$$
This accuracy measure admits a transparent justification: this is the sum, over the players,
of the incentives for a player to change his choice given that other players stick to their choices.
\subsubsection{Special case: saddle points}
An important by its own right particular case of Nash Equilibrium
problem is a {\sl zero sum game}, where $m=2$ and
$\Phi(z)\equiv 0$ (i.e., $\phi_2(z)\equiv-\phi_1(z)$). The convex case of this problem
corresponds to the situation when $\phi(z_1,z_2)\equiv
\phi_1(z_1,z_2)$ is a Lipschitz continuous function which is
convex in $z_1\in Z_1$ and concave in $z_2\in Z_2$, the Nash
equilibria are exactly the saddle points ($\min$ in $z_1$, $\max$
in $z_2$) of $\phi$ on $Z_1\times Z_2$, and the functional
error becomes
$$
\errNash(z_1,z_2)=\max\limits_{(u_1,u_2)\in
Z}\left[\phi(z_1,u_1)-\phi(u_2,z_2)\right].
$$
Recall that the convex-concave
saddle point problem $\min_{z_1\in Z_1}\max_{z_2\in
Z_2}\phi(z_1,z_2)$ gives rise  to the ``primal-dual'' pair
of convex optimization problems
\[
(P): \min_{z_1\in Z_1}\overline{\phi}(z_1),\;\;\;\quad(D): \max_{z_2\in
Z_2}\underline{\phi}(z_2),
\]
where
\[
\overline{\phi}(z_1)=\max_{z_2\in
Z_2}\phi(z_1,z_2),\;\;\;\underline{\phi}(z_2)=\min_{z_1\in Z_1}\phi(z_1,z_2).
\]
The optimal values $\Opt(P)$ and $\Opt(D)$ in these problems are
equal, the set of saddle points of $\phi$ (i.e., the
set of Nash equilibria of the underlying convex Nash problem)  is exactly
the direct product of the optimal sets of $(P)$ and $(D)$, and
$\errNash(z_1,z_2)$ is nothing but the sum of non-optimalities of
$z_1$, $z_2$ considered as approximate solutions to respective
optimization problems:
$$
\errNash(z_1,z_2)=\left[\overline{\phi}(z_1)-\Opt(P)\right]+
\left[\Opt(D)-\underline{\phi}(z_2)\right].
$$

In the sequel, we refer to the v.i. (\ref{VI}) coming from a
convex Nash Equilibrium problem as {\sl Nash} v.i., and to
the just outlined particular case as the {\sl Saddle Point} v.i.
It is easy to verify that in the Saddle Point case
the functional error $\errNash(z)$ is
$\leq \errvi(z)$; this is  not necessary so for a
general Nash v.i.
\subsection{Composite Optimization problem and its saddle point reformulation}\label{scp}
While the algorithm we intend to develop is applicable to a general-type stochastic v.i. with monotone operator,
the applications to be considered in this paper deal with (saddle point reformulation of) {\sl convex composite  optimization problem} (cf. \cite{NemYu}).
\par
As the simplest motivating example, one can keep in mind the minimax problem
\begin{equation}
\label{minmax}
\min\limits_{x\in X} \max\limits_{1\leq i\leq m}\phi_\ell(x),
\end{equation}
where $X\subset\bR^n$ is a convex compact set and $\phi_\ell(x)$ are Lipschitz continuous convex functions on $X$. This problem
can be rewritten as the saddle point problem
\begin{equation}\label{thesameminmax}
\min\limits_{x\in X}\max\limits_{y\in Y}\phi(x,y):=\sum_{ell=1}^my_\ell\phi_\ell(x),
\end{equation}
where $Y=\{y\in\bR^m_+:\sum_{\ell=1}^my_\ell=1\}$ is the standard simplex. The advantages of the saddle point reformulation are twofold. First, when all $\phi_\ell$ are smooth, so is $\phi$, in contrast to the objective in (\ref{minmax}) which typically is nonsmooth; this makes the saddle point reformulation better suited for processing by first order algorithms.  Starting with the breakthrough paper of Nesterov \cite{nes1}, this phenomenon, in its general form, is utilized in the fastest known so far first order algorithms for ``well-structured'' nonsmooth convex programs. Second, in the stochastic case, stochastic oracles providing unbiased estimates of the first order information on $\phi_i$ oracles, while not induce
a similar oracle for the objective of (\ref{minmax}), do induce such an oracle for the v.i. associated with (\ref{thesameminmax}) and thus make the problem amenable to first order algorithms.\par
\subsubsection{Composite minimization problem.}\label{sectcomp} In this paper, we focus on a substantial extension of the minimax problem (\ref{minmax}), namely, on a {\sl Composite minimization problem}
\begin{equation}\label{composite}
\min\limits_{x\in X} \phi(x):=\Phi(\phi_1(x),...,\phi_m(x)),
\end{equation}
where the inner functions $\phi_\ell(\cdot)$ are vector-valued, and the outer function $\Phi$
is real-valued. We are about to impose structural restrictions which allow to reformulate the problem as a ``good'' convex-concave saddle point problem, specifically, as follows:
\begin{enumerate}
\item[{\bf A.}] $X\subset\cX$ is a convex compact;
\item[{\bf B.}] $\phi_\ell(x):X\to \cE_\ell$, $1\leq \ell\leq m$, are Lipschitz
continuous mappings taking values in Euclidean spaces $\cE_\ell$
equipped
with closed convex cones $K_\ell$. We assume $\phi_\ell$ to be
$K_\ell$-convex, meaning that for any $x,x'\in X,\;\lambda\in[0,1]$,
\[
 \phi_\ell(\lambda
x+(1-\lambda)x')\leq_{K_\ell}\lambda\phi_\ell(x)+(1-\lambda)\phi_\ell(x'),
\]
where the notation $a\leq_{K} b\Leftrightarrow b\geq_K a$ means that $b-a\in
K$.
\item[{\bf C.}] $\Phi(\cdot)$ is a convex function on
$\cE=\cE_1\times...\times \cE_m$ given by the Fenchel-type
representation
\begin{equation}\label{representation}
\Phi(u_1,...,u_m)=\max\limits_{y\in Y} \left\{{\sum}_{\ell=1}^m\langle u_\ell,\bA_\ell y+b_\ell\rangle_{\cE_\ell}-\Phi_*(y)\right\},
\end{equation}
for $u_\ell\in \cE_\ell,\;1\leq\ell\leq m$. Here
\\
--- $Y\subset \cY$ is a convex compact
set,\\
--- the affine mappings
$y\mapsto \bA_\ell y+b_\ell:\cY\to \cE_\ell$ are such that
$\bA_\ell y+b_\ell\in K_\ell^*$ for all $y\in Y$ and all $\ell$, $K_\ell^*$ being the cone dual to $K_\ell$,
\\
--- $\Phi_*(y)$ is a given Lipschitz continuous convex function on $Y$.
\end{enumerate}
\par\noindent
Under these assumptions, the
optimization problem (\ref{composite}) is nothing but the primal
problem associated with the saddle point problem
\begin{equation}\label{SP}
\min_{x\in X} \max_{y\in
Y}\left[\phi(x,y)={\sum}_{\ell=1}^m\langle \phi_\ell(x),\bA_\ell y+b_\ell\rangle_{\cE_\ell}-\Phi_*(y)\right]
\end{equation}
and the cost function in the latter problem is Lipschitz continuous and convex-concave due to the convexity of $\Phi_*$,  $K_\ell$-convexity of $\phi_\ell(\cdot)$ and the condition
$\bA_\ell y+b_\ell\in K_\ell^*$ whenever $y\in Y$. The associated Nash v.i.
is given by the domain $Z=X\times Y$ and the monotone mapping
\begin{equation}\label{associated}
F(z)\equiv
F(x,y)=\left[{\sum}_{\ell=1}^m[\phi_\ell^\prime(x)]^*[\bA_\ell y + b_\ell];\;-{\sum}_{\ell=1}^m
\bA_\ell^*\phi_\ell(x)+\Phi_*^\prime(y)\right].
\end{equation}
Same as in the case of minimax problem (\ref{minmax}), the advantage of the saddle point reformulation (\ref{SP}) of (\ref{composite}) is
that, independently of whether $\Phi$ is smooth,  $\phi$ is smooth whenever all $\phi_\ell$ are so. Another advantage, instrumental in the stochastic case, is that $F$ is linear in $\phi_\ell(\cdot)$, so that stochastic oracles providing unbiased estimates of the first order information on $\phi_\ell$ induce straightforwardly an unbiased SO for $F$.
\subsubsection{Example: Matrix Minimax problem}\label{SMMP} For
$1\leq\ell\leq m$, let $\cE_\ell=\bS^{p_\ell}$ be the space of
symmetric $p_\ell\times p_\ell$  matrices equipped with the
Frobenius inner product $\langle A,B\rangle_F=\Tr(AB)$, and let $K_\ell$ be the cone
$\bS^{p_\ell}_+$ of symmetric positive semidefinite $p_\ell\times
p_\ell$ matrices. Now let $X\subset\cX$ be a convex compact set, and $\phi_\ell:X\to\cE_\ell$ be $\bS^{p_\ell}_+$-convex Lipschitz continuous mappings. These data induce the {\sl Matrix Minimax problem}
$$
 \min_{x\in X}\max\limits_{1\leq j\leq
k}\lambda_{\max}\left({\sum}_{\ell=1}^m P_{j\ell}^T \phi_\ell(x) P_{j\ell}\right), \eqno{(P)}
$$ where
$P_{j\ell}$ are given $p_\ell\times q_j$ matrices, and $\lambda_{\max}(A)$ is the maximal eigenvalue of a symmetric matrix $A$.
Observing that for a symmetric $q\times q$ matrix $A$ one has
$$
\lambda_{\max}(A)=\max\limits_{S\in\cS_q}\Tr(AS)
$$
where $\cS_q=\{S\in \bS^q_+: \Tr(S)=1\}$. When
denoting by $Y$ the set of all symmetric positive semidefinite block-di\-a\-go\-nal matrices $y=\Diag\{y_1,...,y_k\}$ with unit trace and diagonal blocks $y_j$ of sizes $q_j\times q_j$, we can represent $(P)$ in the form of
(\ref{composite}), (\ref{representation}) with
$$
\begin{array}{rl}
\Phi(u):=&\max\limits_{1\leq j\leq k}\lambda_{\max}\left({\sum}_{\ell=1}^m
P_{j\ell} u_\ell P_{j\ell}^T\right)
=\max\limits_{y\in Y}{\sum}_{j=1}^k\Tr\left({\sum}_{\ell=1}^mP_{j\ell}^Tu_\ell P_{j\ell}y_j\right)\\
=&\max\limits_{y\in Y} {\sum}_{\ell=1}^m \Tr\left(u_\ell
\left[{\sum}_{j=1}^k P_{j\ell}^Ty_jP_{j\ell}\right]\right)
=\max\limits_{y\in Y}{\sum}_{\ell=1}^m \langle u_\ell,
\bA_\ell y\rangle_F,\\
\bA_{\ell} y=&{\sum}_{j=1}^k P_{j\ell}y_jP_{j\ell}^T\\
\end{array}
$$
Observe that in the simplest case of $k=m$, $p_j=q_j$, $1\leq j\leq m$ and $P_{j\ell}$ equal to $I_p$ for $j=\ell$ and to 0 otherwise,  the problem becomes
\begin{equation}\label{becomesA}
\min_{x\in X}\left[\max_{1\leq \ell\leq
m}\lambda_{\max}(\phi_\ell(x))\right].
\end{equation}
If, in addition, $p_j=q_j=1$ for all $j$, we arrive at the convex minimax problem (\ref{minmax}).
\paragraph{Illustration: Semidefinite Feasibility problem.} With $X$ and $\phi_\ell$ as above, consider the Semidefinite Feasibility problem
$$
\hbox{find\ } x\in X: \psi_\ell(x)\preceq 0,\,1\leq\ell\leq m. \eqno{(S)}
$$
Choosing somehow  scaling factors $\beta_\ell>0$ and setting $\phi_\ell(x)=\beta_\ell\psi_\ell(x)$, we can pose $(\cS)$
as the Matrix Minimax problem
$
\min\limits_{x\in X} \max\limits_{1\leq \ell\leq m}\lambda_{\max}(\phi_\ell(x));
$
$(\cS)$ is solvable if and only if the optimal value in the Matrix Minimax problem is $\leq 0$.
\section{Stochastic Mirror-Prox algorithm}
\label{sect:SMP}
 We are about to present the stochastic version of the deterministic Mirror-Prox algorithm proposed in
\cite{MP}. The method is aimed at solving v.i. (\ref{VI}) associated with a convex compact set $Z\subset\cE$ and a {\sl bounded} monotone operator $F:Z\to\cE$. In contrast to the original version of the method,
below we allow for errors when computing the values of $F$ -- we assume that given a point $z\in Z$, we can
compute an approximation (perhaps random) $\widehat{F}(z)\in \cE$ of $F(z)$.
\subsection{Algorithm's setup}
\paragraph{The setup} for SMP (Stochastic Mirror Prox algorithm) is given by
\begin{enumerate}
\item a norm $\|\cdot\|$ on $\cE$; $\|\cdot\|_*$ stands for
the conjugate norm, see \rf{normstar};
\item a {\sl
distance-generating function} (d.-g.f.) for $Z$, that  is, a
continuous convex function $\omega(\cdot):Z\to\bR$ such that
\begin{enumerate}
\item with $Z^o$ being the set of all points $z\in Z$ such that the
subdifferential $\partial\omega(z)$ of $\omega(\cdot)$ at $z$ is
nonempty, $\partial\omega(\cdot)$ admits a continuous
selection on $Z^o$: there exists a continuous on $Z^o$
vector-valued function $\omega'(z)$ such that
$\omega'(z)\in\partial\omega(z)$ for all $z\in Z^o$;
\item $\omega(\cdot)$ is strongly convex, modulus
$1$, w.r.t. the norm $\|\cdot\|$:
\begin{equation}\label{omega}
\forall (z,z'\in Z^o):
\langle\omega'(z)-\omega'(z'),z-z'\rangle\geq \|z-z'\|^2.
\end{equation}
\end{enumerate}
\end{enumerate}
In order for the SMP associated with the outlined setup to be practical,
$\omega(\cdot)$ and $Z$ should
``fit'' each other, meaning that  one can easily solve problems of
the form
\begin{equation}\label{neweq11}
\min\limits_{z\in Z} \left[\omega(z)+\langle
e,z\rangle\right],\quad e\in \cE.
\end{equation}
\paragraph{The prox-function} associated with a setup for SMP is defined as
\[
V(z,u)=\omega(u)-\omega(z)-\langle\omega'(z),u-z\rangle:Z^o\times
Z\to \bR^+.
\]
We set
\begin{equation}\label{weset}
\begin{array}{lrcl}
(a)&\Theta(z)&=&\max_{u\in Z}V(z,u)\quad [z\in Z^o];\\
(b)&z_{\rm c}&=&\argmin_{Z} \omega(z);\\
(c)&\Omega&=&\sqrt{2\Theta(z_{\rm c})}.\\
\end{array}
\end{equation}
Note that $z_{\rm c}$ is well defined (since $Z$ is a convex compact set
and $\omega(\cdot)$ is continuous and strongly convex on $Z$) and
belongs to $Z^o$ (since $0\in\partial\omega(z_{\rm c})$). Note also that
due to the strong convexity of $\omega$ and the origin of $z_{\rm c}$ we have
\begin{equation}\label{notealso}
\forall (u\in Z):{1\over 2}\|u-z_{\rm c}\|^2\leq \Theta(z_{\rm c})\leq
\max_{z\in Z}\omega(z)-\omega(z_{\rm c});
\end{equation}
in particular we see that
\begin{equation}\label{eqwesee}
Z\subset\{z:\|z-z_{\rm c}\|\leq\Omega\}.
\end{equation}
\paragraph{Prox-mapping.} Given a setup for SMP and a point $z\in Z^o$, we define the associated {\sl prox-mapping} as
$$
Pz(\xi)=\argmin_{u\in Z}\left\{\omega(u)+\langle
\xi-\omega'(z),u\rangle\right\}\equiv\argmin_{u\in
Z}\left\{V(z,u)+\langle \xi,u\rangle\right\}: \cE\to Z^o.
$$
Since $\cS$ is compact and $\omega(\cdot)$ is continuous and strongly convex on $Z$This mapping is clearly well defined.
\subsection{Basic SMP setups}\label{setups}
 We illustrate the just-defined notions with three basic examples.
 \paragraph{Example 1: Euclidean setup.} Here $\cE$ is $\bR^N$ with the standard inner product, $\|\cdot\|_2$ is the
standard Euclidean norm on $\bR^N$ (so that $\|\cdot\|_*=\|\cdot\|$) and $\omega(z)={1\over 2}z^Tz$ (i.e., $Z^o=Z$).
Assume for the sake of simplicity that
$0\in Z$. Then $z_{\rm c}=0$ and $\Omega=\max_{z\in Z}\|z\|_2^2$. The prox-function and the prox-mapping are given by
$
V(z,u)={1\over 2}\|z-u\|_2^2$, $P_z(\xi)=\argmin_{u\in Z} \|(z-\xi)-u\|_2$.
\paragraph{Example 2: Simplex setup.} Here $\cE$ is $\bR^N$, $N>1$, with the standard inner product,
$\|z\|=\|z\|_1:={\sum}_{j=1}^N|z_j|$ (so that $\|\xi\|_*=\max_j|\xi_j|$),
$Z$ is a closed convex subset of the standard simplex $$\cD_N=\{z\in\bR^N:z\geq0,{\sum}_{j=1}^Nz_j=1\}$$
containing its barycenter,
and $\omega(z)={\sum}_{j=1}^Nz_j\ln z_j$ is the entropy. Then
\[
Z^o=\{z\in Z:z>0\}\;\;\mbox{and}\;\;\omega'(z)=[1+\ln z_1;...;1+\ln z_N], \;\;z\in Z^o.
\]
It is easily seen (see, e.g., \cite{Shapiroetal}) that here
\[
z_{\rm c}=[1/N;...;1/N],\;\;\Omega\leq \sqrt{2\ln(N)}\]
(the latter inequality becomes equality when $Z$ contains a vertex of $\cD_N$).
The prox-function is $$V(z,u)={\sum}_{j=1}^Nu_j\ln(u_j/z_j),$$
and the prox-mapping is easy to compute when $Z=\cD_N$:
\[
(P_z(\xi))_j=\left({\sum}_{i=1}^N z_i\exp\{-\xi_i\}\right)^{-1}z_j\exp\{-\xi_j\}.\]

\paragraph{Example 3: Spectahedron setup.} This is the ``matrix analogy'' of the Simplex setup. Specifically, now $\cE$ is the space of
$N\times N$ block-diagonal symmetric matrices, $N>1$, of a given
block-diagonal structure equipped with the Frobenius inner product
$\langle a,b\rangle_F=\Tr(ab)$ and the trace norm
$|a|_1={\sum}_{i=1}^N|\lambda_i(a)|$, where $\lambda_1(a)\geq...\geq
\lambda_N(a)$ are the eigenvalues of a symmetric $N\times N$
matrix $a$; the conjugate norm $|a|_\infty$ is the usual
spectral norm (the largest singular value) of $a$. $Z$ is
assumed to be a closed convex subset of the {\sl spectahedron}
$\cS=\{z\in \cE:z\succeq0,\;\Tr(z)=1\}$ containing the matrix
$N^{-1}I_N$. The d.-g.f. is twice the matrix
entropy $$\omega(z)=2{\sum}_{j=1}^N\lambda_j(z)\ln\lambda_j(z),$$ so that
$Z^o=\{z\in Z: z\succ0\}$ and  $\omega'(z)=2\ln(z)+2I_N$. This setup, similarly to the Simplex
one, results in $z_{\rm c}=N^{-1}I_N$ and $\Omega\leq2\sqrt{\ln
N}$ \cite{NERML}. When $Z=\cS$, it is
relatively easy to compute the prox-mapping (see \cite{NERML,MP});
this task reduces to the singular value decomposition of a matrix
from $\cE$. It should be added that the matrices from $\cS$ are exactly
the matrices of the form
$$a=\cH(b)\equiv(\Tr(\exp\{b\}))^{-1}\exp\{b\}$$
with $b\in \cE$. Note also that when $Z=\cS$, the prox-mapping becomes ``linear in matrix logarithm'': if $z=\cH(a)$, then $P_z(\xi)=\cH(a-\xi/2)$.

\subsection{Algorithm: the construction}
The $t$-step SMP algorithm is applied to the v.i. (\ref{VI}), works as follows:
\begin{algorithm}\label{Basic}
\begin{enumerate}
\item \underline{Initialization:} Choose $r_0\in Z^o$ and
stepsizes $\gamma_\tau>0$, $1\leq\tau\leq t$. \item
\underline{Step $\tau$, $\tau=1,2,...,t$:} Given $r_{\tau-1}\in
Z^o$,  set \be \left\{
\begin{array}{lcl}
w_\tau&=&P_{r_{\tau-1}}(\g_\tau \wh{F}(r_{\tau-1})), \\
r_\tau&=&P_{r_{\tau-1}}(\g_\tau \wh{F}(w_{\tau}))\\
\end{array}
\right. \ee{neq104} When $\tau<t$, loop to step $t+1$. \item At
step $t$, output \begin{equation}\label{output}
\widehat{z}_t=\left[{\sum}_{\tau=1}^t\gamma_\tau\right]^{-1}{\sum}_{\tau=1}^t\gamma_\tau
w_\tau.
\end{equation}
\end{enumerate}
\end{algorithm}
\noindent
Here $\widehat{F}$ is the approximation of $F(\cdot)$ available to the algorithm, so that $\widehat{F}(z)\in\cE$
is the output of the ``black box'' -- the oracle -- representing $F$, the input to the oracle being $z\in Z$. In what follows we
assume that $F$ is a bounded monotone operator represented by a {\sl Stochastic Oracle}.
\paragraph{Stochastic Oracle (SO).}  At the $i$th call to the SO, the input being
$z\in Z$, the oracle returns the vector
$\widehat{F}=\Xi(z,\zeta_i)$, where
$\{\zeta_i\in\bR^N\}_{i=1}^\infty$ is a sequence of i.i.d. random
variables, and $\Xi(z,\zeta):\;Z\times \bR^N\to \cE$ is a Borel
function satisfying the following
\par
{\bf Assumption I}: {\sl With some $\mu\in[0,\infty)$, for all
$z\in Z$  we have}
\begin{equation}\label{eq:expect_w}
\begin{array}{ll}
(a)&\left\|\bE\left\{\Xi(z,\zeta_i)-F(z)\right\}\,\right\|_*\le \mu\\
(b)&\bE\left\{\|\Xi(z,\zeta_i)-F(z)\|_*^2\right\}\leq \sigma^2.\\
\end{array}
\end{equation}
which is slightly milder than \rf{oracle}. The associated version of Algorithm \ref{Basic} will
be referred to as {\sl Stochastic Mirror Prox} (SMP) algorithm.
\par
In some cases, we augment  Assumption I by the following\par {\bf
Assumption II}: {\sl For all $z\in Z$ and all $i$ we have
\begin{equation}\label{eq:expect}
\bE\left\{\exp\{\|\Xi(z,\zeta_i)-F(z)\|_*^2/\sigma^2\}\right\}\leq\exp\{1\}.
\end{equation}
Note that Assumption II implies {\rm
(\ref{eq:expect_w}.$b$)}, since
$$
\exp\{\bE\left\{\|\Xi(z,\zeta_i)-F(z)\|_*^2/\sigma^2\right\}\} \leq
\bE\left\{\exp\{\|\Xi(z,\zeta_i)-F(z)\|_*^2/\sigma^2\}\right\}
$$
by the Jensen inequality.}
\subsection{Algorithm: Main result}

\par
From now on, assume that
the starting point $r_0$ in Algorithm \ref{Basic} is the minimizer
$z_{\rm c}$ of $\omega(\cdot)$ on $Z$. Further, to avoid unnecessarily complicated formulas  (and with no harm to the efficiency estimates)
we stick to the constant stepsize policy $\gamma_\tau\equiv\gamma$,
$1\leq \tau\leq t$, where $t$ is a fixed in advance number of
iterations of the algorithm.
%
 Our main result is as follows:
\begin{theorem}\label{themain} Let v.i. {\rm (\ref{VI})} with monotone operator $F$ satisfying {\rm
\rf{F}} be solved by $t$-step Algorithm \ref{Basic} using a
SO, and let the stepsizes
$\gamma_\tau\equiv \gamma$, $1\leq\tau\leq t$, satisfy
$0<\gamma\leq {1\over \sqrt{3} L}$.
Then
\par
{\rm (i)} Under Assumption I, one has
\begin{equation}\label{K0}
\bE\left\{\errvi(\widehat{z}_t)\right\}\leq K_0(t)\equiv
\left[{\Omega^2\over t\gamma}+{7\gamma\over2}[M^2+2\sigma^2]\right]+
2\mu\Omega,
\end{equation}
where $M$ is the constant from {\rm \rf{F}} and $\Omega$ is given
by {\rm (\ref{weset})}.
\par
{\rm (ii)} Under Assumptions I, II, one has, in addition to {\rm
(\ref{K0})}, for any $\Lambda>0$,
\begin{equation}\label{K1}
\Prob\left\{\errvi(\widehat{z}_t)>K_0(t)+\Lambda
K_1(t)\right\}\leq\exp\{-\Lambda^2/3\}+\exp\{-\Lambda t\},
\end{equation}
where
\[
K_1(t)={7\sigma^2\gamma\over 2}+{2\sigma\Omega \over \sqrt{t}}.
\]
In the case of a Nash v.i., $\errvi(\cdot)$ in {\rm (\ref{K0}),
(\ref{K1})} can be replaced with $\errNash(\cdot)$.
\end{theorem}
\par
When optimizing the bound (\ref{K0}) in $\gamma$, we get
the following
\begin{corollary}\label{maincor} In the situation of
Theorem \ref{themain}, let the stepsizes $\gamma_\tau\equiv \gamma$  be chosen
according to
\begin{equation}\label{optimalstep}
\gamma=\min\left[{1\over \sqrt{3}L},\Omega\sqrt{{2\over 7t(M^2+2\sigma^2)}}\right].
\end{equation}
Then under Assumption I one has
\begin{equation}\label{K0then}
\begin{array}{c}
\bE\left\{\errvi(\widehat{z}_t)\right\}\leq K_0^*(t)\equiv\max
\left[{7\over 4}{\Omega^2 L\over t},\;7\Omega\sqrt{M^2+2\sigma^2\over 3t}\right]+ 2\mu\Omega,\\
\end{array}
\end{equation}
(see {\rm (\ref{weset})}). Under Assumptions I, II, one has, in
addition to {\rm (\ref{K0then})}, for any $\Lambda>0$,
\begin{equation}\label{K1then}
\Prob\left\{\errvi(\widehat{z}_t)>K_0^*(t)+\Lambda
K_1^*(t)\right\}\leq\exp\{-\Lambda^2/3\}+\exp\{-\Lambda t\}
\end{equation}
with
\be
K_1^*(t)={7\over 2}{\Omega \sigma\over\sqrt{t}}.
\ee{K1}
In the case of a Nash v.i., $\errvi(\cdot)$ in {\rm
(\ref{K0then}), (\ref{K1then})} can be replaced with
$\errNash(\cdot)$. \end{corollary}
\begin{remark}
Observe that the upper bound {\rm (\ref{K0then})} for the error of Algorithm \ref{Basic} with stepsize strategy {\rm (\ref{optimalstep})}, in agreement with the lower bound of {\rm \cite{NemYu}}, depends in the same way on the ``size'' $\sigma$ of the perturbation $\Xi(z,\zeta_i)-F(z)$ and on the bound $M$ for the non-Lipschitz component of $F$. From now on to simplify the presentation, with slight abuse of notation, we denote $M$ the maximum of these quantities. Clearly, the latter implies that the bounds {\rm(\ref{eq:expect_w}.b)} and {\rm \rf{expect}}, and thus the bounds {\rm (\ref{K0then}) -- \rf{K1}} of Corollary \ref{maincor} hold with $M$  substituted for $\sigma$.
\end{remark}
\subsection{Comparison with Robust Mirror
SA Algorithm}
Consider the case of a Nash  s.v.i. with
operator $F$ satisfying \rf{F} with $L=0$, and let the
SO be unbiased (i.e., $\mu=0$). In this case, the
bound (\ref{K0then}) reads
\begin{equation}\label{reads}
\bE\left\{\errNash(\widehat{z}_t)\right\}\leq {7\Omega
M\over\sqrt{t}},
\end{equation}
where
$$
M^2=\max\left[\sup_{z,z'\in Z}\|F(z)-F(z')\|_*^2, \;\sup_{z\in
Z}\bE\left\{\|\Xi(z,\zeta_i)-F(z)\|_*^2\right\}\right]
$$
The bound (\ref{reads}) looks very much like the efficiency estimate
\begin{equation}\label{old}
\bE\left\{\errNash(\tilde{z}_t)\right\}\leq O(1){\Omega
\overline{M}\over\sqrt{t}}
\end{equation}
(from now on, all $O(1)$'s are appropriate absolute positive  constants)
for the approximate solution $\tilde{z}_t$ of the $t$-step
{\sl Robust Mirror SA} (RMSA) algorithm
\cite{Shapiroetal}\footnote{$^)$ In this reference, only the
Minimization and the Saddle Point problems are considered. However, the results of \cite{Shapiroetal}
can be easily extended to s.v.i.'s.}$^)$. In the
latter estimate, $\Omega$ is exactly the same as in (\ref{reads}),
and $\overline{M}$ is given by
$$
\overline{M}^2=\max\left[\sup_{z}\|F(z)\|_*^2; \;\sup_{z\in
Z}\bE\left\{\|\Xi(z,\zeta_i)-F(z)\|_*^2\right\}\right].
$$
Note that we always have $M\leq 2\overline{M}$, and typically $M$ and $\overline{M}$
are of the same order of magnitude; it may happen, however (think of the case when $F$ is ``almost constant''), that
$M\ll \overline{M}$. Thus, the bound (\ref{reads}) never is worse, and
sometimes can be much better than the SA bound (\ref{old}). It
should be added that as far as implementation is concerned, the
SMP algorithm is not more complicated than the RMSA (cf. the
description of Algorithm \ref{Basic} with the description
\bse
r_t&=&P_{r_{t-1}}(\widehat{F}(r_{t-1})),\\
\widehat{z}_t&=&\left[{\sum}_{\tau=1}^t\gamma_\tau\right]^{-1}
{\sum}_{\tau=1}^t\gamma_\tau r_\tau,
\ese
 of the RMSA).
\par
The just outlined advantage of SMP as compared to the usual
Stochastic Approximation is not that important, since
``typically'' $M$ and $\overline{M}$ are of the same order. We believe that
the most interesting feature of the SMP algorithm is its ability to take advantage
of a specific structure of a stochastic optimization problem,
namely, insensitivity to the presence in the objective of
large, but smooth and well-observable components.
\par
We are about to consider several less straightforward applications
of the outlined insensitivity of the SMP algorithm to smooth
well-observed components in the objective.
\section{Application to Stochastic Composite minimization}\label{SAa}
Our present goal is to apply the SMP algorithm to the Composite minimization problem (\ref{composite})
in the case when the associated monotone operator (\ref{associated}) is given by a Stochastic Oracle.
\par
{\sl
Throughout this section, the structural assumptions {\bf A} -- {\bf C} from Section \ref{scp} are in force.}
\subsection{Assumptions}\label{assumptions}
We start with augmenting the description of the problem of interest (\ref{composite}), see Section \ref{scp},
with additional assumptions specifying the SO and the SMP setup for the v.i. reformulation
\begin{equation}\label{VIref}
\hbox{find\ } z_*\in Z:=X\times Y: \langle F(z),z-z_*\rangle\geq0\,\,\forall z\in Z
\end{equation}
of (\ref{composite}); here $F$ is the monotone operator (\ref{associated}). Specifically, we assume that
\begin{enumerate}
\item[{\bf D.}] The embedding space $\cX$ of $X$
is equipped with a norm $\|\cdot\|_x$, and $X$ itself -- with a
d.-g.f. $\omega_x(x)$, the associated parameter (\ref{weset}.$c$) being some $\Omega_x$;
\item[{\bf E.}] The spaces $\cE_\ell$, $1\leq \ell\leq m$, where the functions $\phi_\ell$ take their values,
are equipped with norms (not necessarily the Euclidean ones)
$\|\cdot\|_{(\ell)}$ with conjugates $\|\cdot\|_{(\ell,*)}$ such that
\begin{equation}\label{eq880}
\forall v,v'\in X:\left\{
\begin{array}{ll}
(a)&\|[\phi_\ell^\prime(v)-\phi_\ell^\prime(v')]h\|_{(\ell)}\leq [L_x\|v-v'\|_x+M_x]\|h\|_x\\
(b)&\|[\phi_\ell^\prime(v)]h\|_{(\ell)}\leq \Omega_xL_x\|h\|_x\\
\end{array}\right.
\end{equation}
for certain selections $\phi_\ell^\prime(v)\in\partial^{K_\ell}
\phi_\ell(v)$, $v\in X$\footnote{$^)$ For a $K$-convex function
$\phi:X\to \cE$ ($X\subset\cX$ is convex, $K\subset\cE$ is a
closed convex cone) and $x\in X$, the {\sl $K$-subdifferential}
$\partial^K\phi(x)$ is comprised of all linear mappings $h\mapsto
\mathbf{P}h:\cX\to\cE$ such that $\phi(u)\geq_{K} \phi(x)+\mathbf{P}(u-x)$ for
all $u\in X$. When $\phi$ is Lipschitz continuous on $X$,
$\partial^K\phi(x)\neq\emptyset$ for all $x\in X$; if $\phi$ is
differentiable at $x\in\inter X$ (as it is the case almost
everywhere on $\inter X$), one has ${\partial\phi(x)\over \partial
x}\in\partial^K\phi(x)$.}$^)$ and certain nonnegative constants
$L_x,M_x$.
\item[{\bf F.}] Functions $\phi_\ell(\cdot)$ are represented by
an unbiased SO. At the $i$th call to the oracle, $x\in X$ being the
input, the oracle returns vectors $f_\ell(x,\zeta_i)\in \cE_\ell$ and linear mappings
$\bG_\ell(x,\zeta_i)$ from $\cX$ to $\cE_\ell$, $1\leq\ell\leq m$ ($\{\zeta_i\}$ are i.i.d.
random ``oracle noises'') such that for any $x\in X$ and $i=1,2,...$,
\begin{equation}\label{eq900}
\begin{array}{ll}
(a)&\bE\left\{f_\ell(x,\zeta_i)\right\}=\phi_\ell(x),\;1\leq\ell\leq m\\
(b)&\bE\left\{\max\limits_{1\leq\ell\leq
m}\|f_\ell(x,\zeta_i)-\phi_\ell(x)\|_{(\ell)}^2\right\}\leq
M_x^2\Omega_x^2;\\
(c)&\bE\left\{\bG_\ell(x,\zeta_i)\right\}= \phi_\ell^\prime(x),\;1\leq\ell\leq m,\\
(d)&\bE\left\{\max\limits_{{h\in\cX\,\atop\|h\|_x\leq1}}\|[\bG_\ell(x,\zeta_i)-\phi_\ell^\prime(x)]h\|_{(\ell)}^2\right\}\leq
M_x^2,\;1\leq\ell\leq m.\\
\end{array}
\end{equation}
\item[{\bf G.}]
The data participating in the Fenchel-type representation (\ref{representation}) of $\Phi$ are such that
\begin{enumerate} \item the embedding space $\cY$ of $Y$ is equipped with a norm $\|\cdot\|_y$, and $Y$ itself -
--- with a d.-g.f. $\omega_y(y)$, the associated parameter (\ref{weset}.$c$) being some $\Omega_y$;
\item we have $\|y\|_y\leq2\Omega_y$ for all $y\in Y$;\footnote{This requirement can be ensured by  shifting $Y$ to include the origin and the associated shift in $\omega_y(\cdot)$, see (\ref{eqwesee}).}
\item The convex function $\Phi_*(y)$ is given by the precise deterministic first order oracle, and
\begin{equation}\label{suchthat17}
\|\Phi_*^\prime(y)-\Phi_*^\prime(y')\|_{y,*}\leq L_y\|y-y'\|_y+M_y
\end{equation}
for certain selection $\Phi_*^\prime(y)\in\partial\Phi_*(y)$,
$y\in Y$, and some nonnegative $L_y,M_y$.
\end{enumerate}
\end{enumerate}
\paragraph{Stochastic Oracle for (\ref{VIref}).} The assumptions {\bf E} -- {\bf G} induce an unbiased SO for the operator $F$ in (\ref{VIref}), specifically, the oracle
\begin{equation}\label{induces}
\Xi(x,y,\zeta_i)=\left[\sum_{\ell=1}^m
\bG^*_\ell(x,\zeta_i)[\bA_\ell y+b_\ell];\;-\sum_{\ell=1}^m
\bA_\ell^*f_\ell(x,\zeta_i)+\Phi_*^\prime(y)\right],
\end{equation}
and this is the oracle we will use when solving (\ref{composite}) by the SMP algorithm.
\subsection{Setup for the SMP as applied to {\rm
(\ref{VIref}), (\ref{associated}) }}\label{subsetup} In retrospect, the setup for
SMP we are about to present is kind of the best -- resulting in
the best possible efficiency estimate (\ref{K0then}) -- we can
build from the entities participating in the description of the
problem (\ref{composite}) as given by the assumptions {\bf A} -- {\bf G}. Specifically,  we equip the space
$\cE=\cX\times\cY$ with the norm
$$
\|(x,y)\|\equiv \sqrt{\|x\|_x^2/\Omega_x^2+\|y\|_y^2/\Omega_y^2}\quad\left[\Rightarrow\|(\xi,\eta)\|_*=\sqrt{\Omega_x^2\|\xi\|_{x,*}^2+
\Omega_y^2\|\eta\|_{y,*}^2}\right]
$$
and equip $Z=X\times Y$ with the d.-g.f.
$$
\omega(x,y)={1\over\Omega_x^2}\omega_x(x)+{1\over
\Omega_y^2}\omega_y(y)
$$
(it is immediately seen that $\omega(\cdot)$ indeed is a d.-g.f. w.r.t. $X$, $\|\cdot\|$).
The SMP-related properties of our setup are summarized in the
following
\begin{lemma}\label{summarized} Let
\begin{equation}\label{finallyweset}
\cA=\max_{y\in\cY:\|y\|_y\leq1}{\sum}_{\ell=1}^m\|\bA_\ell
y\|_{(\ell,*)},\,\,\cB={\sum}_{\ell=1}^m\|b_\ell\|_{(\ell,*)}.
\end{equation}\par{\rm (i)} The parameter $\Omega$ associated with $\omega(\cdot)$, $\|\cdot\|$, $Z$ by {\rm (\ref{weset})}, is $\leq \sqrt{2}$.
\par{\rm (ii)} One has
\begin{equation}\label{eq900operator}
\forall (z,z'\in Z): \|F(z)-F(z')\|_*\leq L\|z-z'\|+M,
\end{equation}
where
\bse
L&=&4\cA \Omega_x^2\Omega_yL_x+\Omega_x\cB+\Omega_y^2L_y,\\
M&=&\left[3\cA\Omega_y+\cB\right]\Omega_xM_x+\Omega_yM_y.\\
\ese Besides this,
\begin{equation}\label{eq900oracle}
\forall (z\in Z,i): \bE\left\{\Xi(z,\zeta_i)\right\}=F(z);\quad
\bE\left\{\|\Xi(z,\zeta_i)-F(z)\|_*^2\right\}\leq M^2.
\end{equation}
Finally, when relations {\rm (\ref{eq900}.$b$,$d$)} are strengthened to
\begin{equation}\label{eq900str}
\begin{array}{l}
\bE\left\{\exp\left\{\max\limits_{1\leq\ell\leq
m}\|f_\ell(x,\zeta_i)-\phi_\ell(x)\|_{(\ell)}^2/(\Omega_xM)^2\right\}\right\}\leq
\exp\{1\},\\
\bE\left\{\exp\left\{\max\limits_{{h\in \cX,\atop\|h\|_x\leq
1}}\|[\bG_\ell(x)-\phi_\ell^\prime(x)]h\|_{(\ell)}^2/M^2\right\}\right\}\leq
\exp\{1\},\;1\leq\ell\leq m,
\end{array}
\end{equation}
then 
\begin{equation}\label{eq900oraclestr}
\bE\left\{\exp\{\|\Xi(z,\zeta_i)-F(z)\|_*^2/M^2\}\right\}\leq\exp\{1\}.
\end{equation}
\end{lemma}
Combining Lemma \ref{summarized} with Corollary \ref{maincor}, we
get explicit efficiency estimates for the SMP algorithm as applied
to the Stochastic composite minimization problem
(\ref{composite}); these are nothing than estimates (\ref{K0then}) with $\sigma$, $\Omega$ replaced with $M$, $\sqrt{2}$, respectively.
\subsection{Application to Matrix Minimax problem}
Consider Matrix Minimax problem from Section \ref{SMMP}, that is, the problem
\begin{equation}\label{MMX}
 \min_{x\in X}\max\limits_{1\leq j\leq
k}\lambda_{\max}\left({\sum}_{\ell=1}^m P_{j\ell}^T \phi_\ell(x) P_{j\ell}\right),
\end{equation}
where $\phi_\ell(\cdot):X\to \bS^{p_\ell}$ are $\succeq$-convex Lipschitz continuous mappings and $P_{j\ell}\in\bR^{p_\ell\times q_j}$.
As it was shown in Section \ref{SMMP}, (\ref{MMX}) admits saddle point representation
\begin{equation}\label{SPMMX}
\begin{array}{c}
\min\limits_{x\in X} \max\limits_{y=\sDiag\{y_1,...,y_k\}\in Y}\sum_{\ell=1}^m\langle\phi_\ell(x),\bA_\ell y\rangle_F,\\
\bA_\ell y=\sum_{j=1}^k P_{j\ell}y_jP_{j\ell}^T,\\
\end{array}
\end{equation}
where $Y$ is the spectahedron  in the space $\cY$ of block-diagonal symmetric matrices $y=\Diag\{y_1,...,y_k\}$ with $k$ diagonal blocks
of sizes $q_1,...,q_k$. Note that we are in the situation described by assumptions {\bf A} -- {\bf C} from Section \ref{scp}, with $\Phi_*(\cdot)\equiv 0$, and $K_\ell:=\bS^{p_\ell}_+\subset \cE_\ell:=\bS^{p_\ell}$. Using the spectahedron setup for $Y$  and setting  $M_y=L_y=0$, $\Omega_y=2\sqrt{\ln(q_1+...+q_k)}$, we meet all assumptions in {\bf G}. Let us specify the norms $\|\cdot\|_{(\ell)}$ on the spaces $\cE_\ell=\bS^{p_\ell}$ as the standard matrix norms $|\cdot|_\infty$ (maximal singular value), and let assumptions {\bf D} -- {\bf F} from Section \ref{assumptions} take place. Note that in the case in question  (\ref{finallyweset}) reads
$$
\cA=\max\limits_{1\leq j\leq k}\max\limits_{\xi\in\bR^{q_j}:\|\xi\|_2=1}\sum_{\ell=1}^m\|P_{j\ell}\xi\|_2^2=
\max\limits_{1\leq j\leq k}|\sum_{\ell=1}^mP_{j\ell}^TP_{j\ell}|_\infty,\quad \cB=0.
$$
(look what are the extreme points of $Y$), so that the quantities $L$, $M$ as given by Lemma \ref{summarized} become
\begin{equation}\label{aneq123becomes}
\begin{array}{rcl}
L&=&O(1)\cA\sqrt{\ln(q_1+...+q_k)}\Omega_x^2L_x,\\
M&=&O(1)\cA\sqrt{\ln(q_1+...+q_k)}\Omega_xM_x.\\
\end{array}
\end{equation}
Note that in the case of problem (\ref{becomesA}) one has $\cA=1$.
\subsubsection{Application to Stochastic Semidefinite Feasibility problem}
Now consider the stochastic version of the Semidefinite Feasibility problem $(S)$:
\begin{equation}\label{SDPFeas}
\hbox{find\ } x\in X: \psi_\ell(x)\preceq0,\,1\leq\ell\leq m.
\end{equation}
Assuming form now on that the latter problem is feasible, we can rewrite it as the Matrix Minimax problem
\begin{equation}\label{MatrMin1}
\Opt=\min\limits_{x\in X}\max\limits_{1\leq \ell \leq m} \lambda_{\max}(\phi_\ell(x)),\quad \phi_\ell(x)=\beta_\ell\psi_\ell(x),
\end{equation}
where $\beta_\ell>0$ are ``scale factors'' we are free to choose.
We are about to show how to use this freedom in order to improve the SMP efficiency estimates.
\par
Assuming, same as in the case of a general-type Matrix Minimax problem, that {\bf D} takes place, let us modify assumptions
{\bf E}, {\bf F} as follows:
\par
{\bf E$'$:} the $\succeq$-convex Lipschitz functions $\psi_\ell:X\to\bS^{p_\ell}$ are such that
\begin{equation}\label{suchthat1}
\begin{array}{l}
\max\limits_{h\in\cX,\,\|h\|_x\leq1}|[\psi_\ell^\prime(x)-\psi_\ell^\prime(x')]h|_\infty\leq
L_\ell\|x-x'\|_x+M_\ell,\\
\max\limits_{{h\in\cX,\,\|h\|_x\leq1}}|\psi_\ell^\prime(x)h|_\infty\leq
\Omega_xL_\ell\\
\end{array}
\end{equation}
for certain selections $\psi_\ell^\prime(x)\in\partial^{K_\ell}
\psi_\ell(x)$, $x\in X$,  with some known nonnegative constants
$L_\ell$, $M_\ell$.
\par
{\bf  F$'$:} $\psi_\ell(\cdot)$ are
represented by an SO which at the $i$th call, the input being $x\in
X$, returns the matrices $\widehat{f}_\ell(x,\zeta_i)\in \bS^{p_\ell}$
and the linear maps $\widehat{\bG}_\ell(x,\zeta_i)$ from $\cX$ to
$\bS^{p_\ell}$ ($\{\zeta_i\}$ are i.i.d.
random ``oracle noises'') such that for any $x\in X$ it holds
\begin{equation}\label{suchthat2}
\begin{array}{ll}
(a)&
\bE\left\{\widehat{f}_\ell(x,\zeta_i)\right\}=\psi_\ell(x),\;\;\bE\left\{\widehat{\bG}_\ell(x,\zeta_i)\right\}=\psi_\ell^\prime(x),
\,1\leq\ell\leq
m\\
(b)& \bE\left\{\max\limits_{1\leq\ell\leq m}
|\widehat{f}_\ell(x,\zeta_i)-\psi_\ell(x)|_\infty^2/(\Omega_xM_\ell)^2\right\}\leq
1\\
(c)&\bE\left\{\max\limits_{{h\in\cX,\atop\|h\|_x\leq1}}|[\widehat{\bG}_\ell(x,\zeta_i)-\psi_\ell^\prime(x)]h|_\infty^2/M_\ell^2\right\}\leq
1,\, 1\leq\ell\leq m. \\
\end{array}
\end{equation}
\par
Given a number $t$ of steps of the SMP algorithm, let us act as
follows.
\par {\bf (I):} We compute the $m$ quantities
$\mu_\ell={\Omega_xL_\ell\over\sqrt{t}}+M_\ell$, $\ell=1,...,m$,
and set
\begin{equation}\label{andset} \mu=\max_{1\leq \ell\leq
m}\mu_\ell, \;\;\beta_\ell={\mu\over
\mu_\ell},\phi_\ell(\cdot)=\beta_\ell \psi_\ell(\cdot),\;\;
L_x=\Omega_x^{-1}\mu\sqrt{t},\;\; M_x=\mu. \end{equation}
Note that
by construction $\beta_\ell\geq1$ and  $L_x/L_\ell\geq
\beta_\ell$, $M_x/M_\ell\geq\beta_\ell$ for all $\ell$, so that
the functions $\phi_\ell$ satisfy (\ref{eq880}) with the just
defined $L_x$, $M_x$. Further, the SO for $\psi_\ell(\cdot)$'s can
be converted into an SO for $\phi_\ell(\cdot)$'s by setting
\[
f_\ell(x,\zeta)=\beta_\ell \widehat{f}_\ell(x,\zeta),\;\;
\bG_\ell(x,\zeta)=\beta_\ell \widehat{\bG}_\ell(x,\zeta).
\] By
(\ref{suchthat2}) and due to $L_x/L_\ell\geq\beta_\ell$, $M_x/M_\ell\geq\beta_\ell$, this oracle satisfies (\ref{eq900}).
\par
{\bf (II)} We then build the Stochastic Matrix Minimax problem
\begin{equation}\label{smmp}
\Opt=\min_{x\in X}\max_{1\leq\ell\leq m}\lambda_{\max}(\phi_\ell(x)),
\end{equation}
associated with the just defined $\phi_1,...,\phi_m$ and solve this Stochastic composite problem by $t$-step SMP algorithm.
Combining Lemma \ref{summarized}, Corollary \ref{maincor} and
taking into account the origin of the quantities $L_x$, $M_x$, and the fact that $\cA=1$, $\cB=0$, we
arrive at the following result:
\begin{proposition}\label{wearriveat} With the outlined construction,  the $t$-step SMP algorithm with the setup presented in Section \ref{subsetup} (where one uses $\cA=1,\cB=L_y=M_y=0$ and the just defined $L_x,M_x$) and constant stepsizes
$\gamma_\tau\equiv\gamma$ defined by {\rm (\ref{optimalstep})}, yields an approximate solution
$\widehat{z}_t=(\widehat{x}_t,\widehat{y}_t)$ such that
\begin{equation}\label{accuracy16}
\begin{array}{rcl}
\bE\left\{\max\limits_{1\leq \ell\le m}\max[\beta_\ell\lambda_{\max}(\psi_\ell(\widehat{x}_t),0]\right\}&\leq &\bE\left\{\max\limits_{1\leq\ell\leq m}\beta_\ell
\lambda_{\max}(\psi_\ell(\widehat{x}_t))-\Opt\right\}\\
&\leq& K_0(t)\equiv 80{\Omega_x\mu\sqrt{\ln({\sum}_{\ell=1}^mp_\ell)}\over\sqrt{t}},\\
\end{array}
\end{equation}
(cf. {\rm (\ref{K0then})} and take into account that we are in the
case of $\Omega=\sqrt{2}$, while the optimal value in {\rm
(\ref{smmp})} is nonpositive, since {\rm (\ref{SDPFeas})} is
feasible).
\par
Furthermore, if assumptions {\rm (\ref{suchthat2}.$b$,$c$)} are strengthened to
\bse
\bE\left\{\max\limits_{1\leq\ell\leq m} \exp\{
|\widehat{f}_\ell(x,\zeta_i)-\psi_\ell(x)|_\infty^2/(\Omega_xM_\ell)^2\}\right\}&\leq&
\exp\{1\},\\
\bE\left\{\exp\{\max\limits_{{h\in \cX,\,
\|h\|_x\leq1}}|[\widehat{\bG}_\ell(x,\zeta_i)-\psi_\ell^\prime(x)]h|_\infty^2/M_\ell^2\}\right\}&\leq&
\exp\{1\},\;\; 1\leq\ell\leq m,
\ese
then, in addition to {\rm (\ref{accuracy16})}, we have for any $\Lambda> 0$:
$$
\begin{array}{l}
\Prob\left\{\max\limits_{1\leq\ell\leq
m}\max[\beta_\ell \lambda_{\max}(\psi_\ell(\widehat{x}_t)),0]>K_0(t)+\Lambda K_1(t)\right\}\\
\multicolumn{1}{r}{\leq\exp\{-\Lambda^2/3\}+\exp\{-\Lambda t\},}
\end{array}
$$
where
\[
K_1(t)={15\Omega_x\mu\sqrt{\ln({\sum}_{\ell=1}^mp_\ell)}\over\sqrt{t}}.
\]
\end{proposition}
\paragraph{Discussion} Imagine that instead of solving the system of matrix inequalities (\ref{SDPFeas}), we were interested to solve just a single matrix inequality $\psi_\ell(x)
\preceq 0$, $x\in X$. When solving this inequality by the SMP
algorithm as explained above, the efficiency estimate would be
\bse
\bE\left\{\max[\lambda_{\max}(\psi_\ell(\widehat{x}_t^\ell)),0]\right\}&\leq&
O(1)\sqrt{\ln(p_\ell+1)}\Omega_x\left[{\Omega_xL_\ell\over
t}+{M_\ell\over \sqrt{t}}\right]\\
&=&O(1)\sqrt{\ln(p_\ell+1)}\beta_\ell^{-1}{\Omega_x\mu\over\sqrt{t}},\\
\ese
(recall that the matrix inequality in question is feasible), where
$\widehat{x}^\ell_t$ is the resulting approximate solution.
Looking at (\ref{accuracy16}), we see that the expected accuracy of the SMP as applied, in
the aforementioned manner, to (\ref{SDPFeas}) is only by a logarithmic in ${\sum}_\ell p_\ell$ factor worse:
\begin{equation}\label{logworse}
\begin{array}{rcl}
\bE\left\{\max[\lambda_{\max}(\psi_\ell(\widehat{x}_t),0]\right\}&\leq&
O(1)\sqrt{\ln({\sum}_{\ell=1}^mp_\ell)}\beta_\ell^{-1}{\Omega_x\mu\over\sqrt{t}}\\
&=&O(1)\sqrt{\ln({\sum}_{\ell=1}^mp_\ell)}{\Omega_x\mu_\ell\over\sqrt{t}}.\\
\end{array}
\end{equation}
Thus, as far as the quality of the SPM-generated solution is
concerned, passing from solving a single matrix inequality to
solving a system of $m$ inequalities is  ``nearly costless''. As an illustration, consider the case where
some of $\psi_\ell$ are ``easy'' -- smooth and easy-to-observe
($M_\ell=0$), while the remaining $\psi_\ell$ are ``difficult'', i.e., might
be non-smooth and/or difficult-to-observe ($\Omega_xL_\ell/\sqrt{t}\leq M_\ell$). In this
case, (\ref{logworse}) reads
$$
\bE\left\{\psi_\ell(\widehat{x}_t)\right\}\leq
O(1)\sqrt{\ln({\sum}_{\ell=1}^mp_\ell)}\cdot\left\{
\begin{array}{ll}
{\Omega_x^2 L_\ell\over t},&\psi_\ell\hbox{\ is easy},\\
{\Omega_xM_\ell\over\sqrt{t}},&\psi_\ell\hbox{\ is difficult}.\\
\end{array}
\right.
$$
In other words, the violations of the easy and the difficult constraints
in (\ref{SDPFeas}) converge to 0 as $t\to\infty$ with the rates $O(1/t)$
and $O(1/\sqrt{t})$, respectively. It should be added that when
$X$ is the unit Euclidean ball in $\cX=\bR^n$ and $X,\,\cX$ are
equipped with the Euclidean setup, the rates of
convergence $O(1/t)$ and $O(1/\sqrt{t})$ are the best rates one can
achieve without imposing bounds on $n$ and/or imposing additional
restrictions on $\psi_\ell$'s.
\subsection{Eigenvalue optimization via SMP}
The problem we are interested in now is
\begin{equation}\label{interestednow}
\begin{array}{rcl}
\Opt&=&\min\limits_{x\in X}f(x):=\lambda_{\max}(A_0+x_1A_1+...+x_nA_n),\\
X&=&\{x\in\bR^n:x\geq0,{\sum}_{i=1}^nx_i=1\},
\end{array}
\end{equation}
where $A_0,A_1,...,A_n$, $n>1$, belong to the space $\bS$ of symmetric matrices with block-diagonal structure $(p_1,...,p_m)$ (i.e.,  a matrix $A\in\bS$ is block-diagonal with $p_\ell\times p_\ell$ diagonal blocks $A^\ell$, $1\leq \ell\leq m$). We set
$$
p^{(\kappa)}={\sum}_{\ell=1}^mp_\ell^\kappa,\;\;\kappa=1,2,...; \;\; p^{\max}=\max_\ell p_\ell;\;\;A_\infty=\max\limits_{1\leq j\leq n}|A_j|_\infty.
$$
Setting
$$
\phi_\ell:\, X\mapsto \cE_\ell=\bS^{p_\ell},\;\;\phi_\ell(x)=A_0^\ell+{\sum}_{j=1}^nx_jA_j^\ell,\;\;1\leq\ell\leq m,
$$
we represent (\ref{interestednow}) as a particular case of the Matrix Minimax problem (\ref{becomesA}), with all functions $\phi_\ell(x)$ being affine and $X$ being the standard simplex in $\cX=\bR^n$.
\par
Now, since $A_j$ are known in advance, there is nothing stochastic in our problem, and it can be solved either by interior point methods, or by ``computationally cheap'' gradient-type methods which are preferable when the problem is large-scale and medium accuracy solutions are sought. For instance, one can apply the $t$-step (deterministic) Mirror Prox algorithm DMP from \cite{MP} to the saddle point reformulation (\ref{SP}) of our specific Matrix Minimax problem, i.e., to the saddle point problem
\begin{equation}\label{SPbecomes}\begin{array}{c}
\min\limits_{x\in X}\max\limits_{y\in Y}\,\langle y,A_0+{\sum}_{j=1}^nx_jA_j\rangle_F, \\
Y=\left\{y=\Diag\{y_1,...,y_m\}: y_\ell\in\bS^{p_\ell}_+,\,1\leq\ell\leq m,\,\Tr(y)=1\right\}.
\end{array}
\end{equation}
The accuracy of the approximate solution $\tilde{x}_t$ of the DMP algorithm is \cite[Example 2]{MP}
\begin{equation}\label{effreads}
f(\tilde{x}_t)-\Opt\leq O(1){\sqrt{\ln(n)\ln(p^{(1)})}A_\infty\over t}.
\end{equation}
This efficiency estimate is the best known so far among those attainable with  ``computationally cheap'' deterministic methods. On the other hand, the complexity of one step of the algorithm is dominated, up to an absolute constant factor, by the necessity, given $x\in X$ and $y\in Y$,
\begin{enumerate}
\item to compute $A_0+{\sum}_{j=1}^nx_jA_j$ and $[\Tr(YA_1);...;\Tr(YA_n)]$;
\item to compute the eigenvalue decomposition of an $y\in\bS$.
\end{enumerate}
When using the standard Linear Algebra, the computational effort per step is
\begin{equation}\label{effort}
\cC_{\rm det}=O(1)[n p^{(2)} + p^{(3)}]
\end{equation}
arithmetic operations.
\par
We are about to demonstrate that one can equip the deterministic problem in question by an ``artificial'' SO in such a way that the associated SMP algorithm, under certain circumstances, exhibits better performance than deterministic algorithms. Let us consider the following construction of the SO for $F$ (different from the SO (\ref{induces})!).
Observe that the monotone operator associated with the saddle point problem (\ref{SPbecomes}) is
\begin{equation}\label{Fis}
F(x,y)=\bigg[\underbrace{[\Tr(yA_1);...;\Tr(yA_n)]}_{F^x(x,y)};\;
\underbrace{-A_0-{\sum}_{j=1}^nx_jA_j}_{F^y(x,y)}\bigg].
\end{equation}
Given $x\in X$, $y=\Diag\{y_1,...,y_m\}\in Y$, we build a random estimate $\Xi=[\Xi^x;\Xi^y]$ of $F(x,y)=[F^x(x,y);F^y(x,y)]$ as follows:
\begin{enumerate}
\item we generate a realization $\jmath$ of a random variable taking values $1,...,n$ with probabilities $x_1,...,x_n$ (recall that $x\in X$, the standard simplex, so that $x$ indeed can be seen as a probability distribution), and set
\begin{equation}\label{andset17}
\Xi^y=A_0+A_\jmath;
\end{equation}
\item we compute the quantities $\nu_\ell=\Tr(y_\ell)$, $1\leq\ell\leq m$. Since $y\in Y$, we have $\nu_\ell\geq0$ and ${\sum}_{\ell=1}^m\nu_\ell=1$. We further generate a realization $\imath$ of random variable taking values $1,...,m$ with probabilities $\nu_1,...,\nu_m$, and set
\begin{equation}\label{andset18}
\Xi^x=[\Tr(A_1^\imath\bar{y}_\imath);...;\Tr(A_n^\imath\bar{y}_\imath)],\,\,\bar{y}_\imath=(\Tr(y_\imath))^{-1}y_\imath.
\end{equation}
\end{enumerate}
The just defined random estimate $\Xi$ of $F(x,y)$ can be expressed as a deterministic function $\Xi(x,y,\eta)$ of $(x,y)$
and random variable $\eta$ uniformly distributed on $[0,1]$. Assuming all matrices $A_j$ directly available (so that it takes $O(1)$ arithmetic operations to extract a particular entry of $A_j$ given $j$ and indexes of the entry) and given $x,y$ and $\eta$, the value $\Xi(x,y,\xi)$ can be computed with the arithmetic cost
$O(1)(n(p^{\max})^2+p^{(2)})
$ (indeed, $O(1)(n+p^{(1)})$ operations are needed to convert $\eta$ into $\imath$ and $\jmath$, $O(1)p^{(2)}$ operations are used to write down the $y$-component $-A_0-A_\jmath$ of $\Xi$, and $O(1)n(p^{\max})^2$ operations are needed to compute $\Xi^x$).
Now consider the SO's $\Xi_k$ ($k$ is a positive integer) obtained by averaging the outputs of $k$ calls to our basic oracle $\Xi$. Specifically, at the $i$th call to the oracle $\Xi_k$, $z=(x,y)\in Z=X\times Y$ being the input, the oracle returns the vector
$$
\Xi_k(z,\zeta_i)={1\over k}{\sum}_{s=1}^k \Xi(z,\eta_{is}),
$$
where $\zeta_i=[\eta_{i1};...;\eta_{ik}]$ and $\{\eta_{is}\}_{{1\leq i,\,
1\leq s\leq k}}$ are independent random variables uniformly distributed on $[0,1]$. Note that the arithmetic cost of a single call to $\Xi_k$ is
$$
\cC_k=O(1)k(n(p^{\max})^2+p^{(2)}).
$$
The Nash v.i. associated with (\ref{SPbecomes}) and the stochastic oracle $\Xi_k$ ($k$ is the first parameter of our construction)
specify a Nash s.v.i. on the domain $Z=X\times Y$. We equip the standard simplex $X$ and its embedding space $\cX=\bR^n$ with the Simplex setup, and the spectahedron $Y$ and its embedding space $\bS$ with the Spectahedron setup (see Section \ref{setups}). Let us next combine the $x$- and the $y$-setups, exactly as explained in the beginning of Section \ref{subsetup}, into an SMP setup  for the domain $Z=X\times Y$ -- a {d.-g.f.} $\omega(\cdot)$ and a norm $\|\cdot\|$ on the embedding space $\bR^n\times (\bS^{p_1}\times...\times\bS^{p_\ell})$ of $Z$. The SMP-related properties of the resulting setup are summarized in the following statement.
\begin{lemma}\label{summarized1} Let $n\geq3$, $p^{(1)}\geq3$. Then\par
{\rm (i)} The parameter of the just defined d.-g.f. $\omega$
w.r.t. the just defined norm $\|\cdot\|$ is $\Omega=\sqrt{2}$.
\par{\rm (ii)} For any $z,z'\in Z$ one has
\begin{equation}\label{neweq900operator}
\|F(z)-F(z')\|_*\leq L\|z-z'\|,\;\;\;L= \sqrt{2}\left[\ln(n)+\ln(p^{(1)})\right]A_\infty.
\end{equation}
Besides this, for any $(z\in Z,\;i=1,2,...$,
\begin{equation}\label{neweq900oracle}
\begin{array}{ll}
(a)&\bE\left\{\Xi_k(z,\zeta_i)\right\}=F(z);\\
(b)&
\bE\left\{\exp\{\|\Xi(z,\zeta_i)-F(z)\|_*^2/M^2\}\right\}\leq\exp\{1\},\\
&M=27[\ln(n)+\ln(p^{(1)})]A_\infty/\sqrt{k}.
\\
\end{array}
\end{equation}
\end{lemma}
\par
Combining Lemma \ref{summarized1} and Corollary \ref{maincor}, we arrive at the following
\begin{proposition}\label{thefinal} With properly chosen positive absolute constants $O(1)$, the $t$-step SMP algorithm with constant stepsizes
$$
\gamma_\tau=O(1){\min[1,\sqrt{k/t}]\over \ln(np^{(1)}) A_\infty},\,1\leq\tau\leq t\eqno{[A_\infty=\max\limits_{1\leq j\leq n}|A_j|_\infty]}
$$
as applied to the saddle point reformulation of problem {\rm (\ref{interestednow})}, the stochastic oracle being $\Xi_k$, produces a random feasible approximate solution $\widehat{x}_t$ to the problem with the error
$$
\epsilon(\widehat{x}_t)=\lambda_{\max}\left(A_0+{\sum}_{j=1}^n[\widehat{x}_t]_jA_j\right)-\Opt
$$
satisfying
\begin{equation}\label{eq123654}
\bE \left\{\epsilon(\widehat{x}_t)\right\}\leq O(1)\ln(np^{(1)})A_\infty\left[{1\over t}+{1\over\sqrt{kt}}\right],
\end{equation}
and for any $\Lambda>0$:
\bse
\Prob\left\{\epsilon(\widehat{x}_t)>O(1)\ln(np^{(1)})A_\infty\left[{1\over t}+{1+\Lambda\over\sqrt{kt}}\right]\right\}
\leq \exp\{-\Lambda^2/3\}+\exp\{-\Lambda t\}.
\ese
Further, assuming that all matrices $A_j$ are directly available, the overall computational effort to compute $\widehat{x}_t$ is
\be
\cC=O(1)t\left[k(n(p^{\max})^2+p^{(2)})+p^{(3)}\right]
\ee{compex}
arithmetic operations.
\end{proposition}
\par
To justify the bound \rf{compex} it suffices to note that  $O(1)k(n(p^{\max})^2+p^{(2)})$ operations per step is the price of two  calls to the stochastic oracle $\Xi_k$ and $O(1)(n+p^{(3)})$ operations per step is the price of computing two prox mappings.
\paragraph{Discussion} Let us find out whether randomization can help when solving a large-scale problem (\ref{interestednow}), that is, whether, given quality of the resulting approximate solution, the computational effort to build such a solution with the Stochastic Mirror Prox algorithm SMP can be essentially less than the one for the deterministic Mirror Prox algorithm DMP. To simplify our considerations, assume from now on that $p_\ell=p$, $1\leq \ell\leq m$, and that $\ln(n)=O(1)\ln(mp)$. Assume also that we are interested
in a (perhaps, random) solution $\widehat{x}_t$ which with probability $\geq 1-\delta$ satisfies $\epsilon(\widehat{x}_t)\leq\epsilon$. We fix a tolerance $\delta\ll 1$ and the {\sl relative accuracy} $\nu={\epsilon\over \ln(mnp)A_\infty}\leq1$ and look what happens when (some of) the sizes $m,n,p$ of the problem become large.
\par
Observe first of all that the overall computational effort to solve (\ref{interestednow}) within relative accuracy $\nu$ with the DMP algorithm is
$$
\cC^{\rm DMP}(\nu)=O(1)m(n+p)p^2\nu^{-1}
$$
operations (see (\ref{effreads}), (\ref{effort})). As for the SMP algorithm, let us choose $k$ which balances the per step computational effort $O(1)k(n(p^{\max})^2+p^{(2)})=O(1)k(n+m)p^2$ to produce the answers of the stochastic oracle and the per step cost of prox mappings $O(1)(n+mp^{3})$, that is, let us set $k=\hbox{Ceil}\left({mp\over m+n}\right)$ \footnote{The rationale behind balancing is clear: with the just defined  $k$, the arithmetic cost of an iteration still is of the same order as when $k=1$, while the left hand side in the efficiency estimate (\ref{eq123654})
becomes better than for $k=1$.}. With this choice of $k$,
Proposition \ref{thefinal} says that to get a solution of the required quality, it suffices to carry out
\begin{equation}\label{suffices}
t=O(1)\left[\nu^{-1}+\ln(1/\delta)k^{-1}\nu^{-2}\right]
\end{equation} steps of the method, provided that this number of steps is $\geq\sqrt{\ln(2/\delta)}$. The latter assumption is automatically satisfied when the absolute constant factor in (\ref{suffices}) is $\geq1$ and $\nu\sqrt{\ln(2/\delta)}\leq1$, which we assume from now on. Combining (\ref{suffices}) and the upper bounds on the arithmetic cost of an SMP step stated in Proposition \ref{thefinal}, we conclude that the overall computational effort to produce a solution of the required quality with the SMP algorithm is
$$
\cC^{\rm SMP}(\nu,\delta)=O(1)k(n+m)p^2\left[\nu^{-1}+\ln(1/\delta)k^{-1}\nu^{-2}\right]
$$
operations, so that
$$
R:={\cC^{\rm DMP}(\nu)\over \cC^{\rm SMP}(\nu,\delta)}=O(1){m(n+p)\over (m+n)(k+\ln(1/\delta)\nu^{-1})}\eqno{[k=\hbox{Ceil}\left({mp\over m+n}\right)]}
$$
We see that {\sl when $\nu$, $\delta$ are fixed, $m\geq n/p$ and $n/p$ is large, then $R$ is large as well},
that is, the randomized algorithm  significantly outperforms its deterministic counterpart.
\par
Another interesting observation is as follows.  In order to produce, with probability $\geq 1-\delta$, an approximate solution to (\ref{interestednow}) with relative accuracy $\nu$, the just defined SMP algorithm requires $t$ steps, with $t$ given by (\ref{suffices}), and at every one of these steps it ``visits'' $O(1)k(m+n)p^2$ randomly chosen entries in the data matrices $A_0,A_1,...,A_n$. The overall number of data entries visited by the algorithm is therefore $N^{\rm SMP}=O(1)tk(m+n)p^2=
O(1)\left[k\nu^{-1}+\ln(1/\delta)\nu^{-2}\right](m+n)p^2$. At the same time, the total number of data entries is $N^{\rm tot}=m(n+1)p^2$. Therefore
$$
\vartheta:={N^{\rm SMP}\over N^{\rm tot}}=O(1)\left[k\nu^{-1}+\ln(1/\delta)\nu^{-2}\right]\left[{1\over m}+{1\over n}\right]\eqno{[k=\hbox{Ceil}\left({mp\over m+n}\right)]}
$$
We see that {\sl when $\delta$, $\nu$ are fixed, $m\geq n/p$ and $n/p$ is large,   $\vartheta$ is small, i.e., the approximate solution of the required quality is built when inspecting a tiny fraction of the data.} This {\sl sublinear time behavior} \cite{ronit} was already observed in \cite{Shapiroetal}  for the Robust Mirror Descent Stochastic Approximation as applied to a matrix game (the latter problem is the particular case of (\ref{interestednow}) with $p_1=...=p_m=1$ and $A_0=0$). Note also that an ``ad hoc'' sublinear time algorithm for a matrix game, in retrospect close to the one from \cite{Shapiroetal}, was discovered in \cite{Khach} as early as in 1995.

\section{Appendix}
\subsection{Preliminaries}
We need the following technical result about the algorithm \rf{neq104}, (\ref{output}):
\begin{theorem}\label{prop1prop1} Consider $t$-step algorithm \ref{Basic} as applied to a v.i. {\rm (\ref{VI})}
with a monotone operator $F$ satisfying \rf{F}.  For
$\tau=1,2,...$, let us set
$$\Delta_\tau=F(w_\tau)-\wh{F}(w_\tau);$$
for $z$ belonging to the trajectory $\{r_0, \,w_1, \,r_1, ...,$
$\,w_t, \,r_t\}$ of the algorithm, let $$\e_z=
\|\widehat{F}(z)-F(z)\|_*,$$ and let $\{y_\tau\in Z^o\}_{\tau=0}^t$
be the sequence given by the recurrence
\begin{equation}\label{recurrence}
y_\tau=P_{y_{\tau-1}}(\g_\tau\Delta_\tau),\;\;y_0=r_0.
\end{equation}
Assume that \be \gamma_\tau\leq {1\over \sqrt{3}L},
\ee{gamma}
 Then \be \errvi(\widehat{z}_t) \le
\left({\sum}_{\tau=1}^t\gamma_\tau\right)^{-1}\Gamma(t),
\ee{zneq200} where $\errvi(\widehat{z}_t)$ is defined in {\rm
\rf{error}},
\begin{eqnarray}\label{eq:Gt0}
\Gamma(t)&=& 2\Theta(r_0)+{\sum}_{\tau=1}^t{3\g_\tau^2\over
2}\left[M^2+(\e_{r_{\tau-1}}+\e_{w_\tau})^2+{\e^2_{w_\tau}\over3}\right]\\
&&+{\sum}_{\tau=1}^t \langle
\g_\tau\Delta_\tau,w_\tau-y_{\tau-1}\rangle\nonumber
 \end{eqnarray}
 and $\Theta(\cdot)$ is defined by {\rm (\ref{weset})}.
\par
Finally, when {\rm (\ref{VI})} is a Nash v.i., one
can replace $\errvi(\widehat{z}_t)$ in {\rm \rf{zneq200}} with
$\errNash(\widehat{z}_t)$.
\end{theorem}
\paragraph{Proof of Theorem \ref{prop1prop1}}
\par
{\bf 1$^0$.}
We start with the following simple observation: if $r_e$ is a
solution to {\rm (\ref{neweq11})}, then $\partial_Z \omega(r_e)$
contains $-e$ and thus is nonempty, so that $r_e\in Z^o$.
Moreover, one has
\begin{equation}\label{neweq12}
\langle \omega'(r_e)+e,u-r_e\rangle \geq0\,\,\forall u\in Z.
\end{equation}
Indeed, by continuity argument, it suffices to verify the
inequality in the case when $u\in\rint (Z)\subset Z^o$. For such
an $u$, the convex function \[ f(t)=\omega(r_e+t(u-r_e))+\langle
r_e+t(u-r_e),e\rangle,\;\; t\in[0,1]\]
 is continuous on $[0,1]$ and
has a continuous on $[0,1]$ field of subgradients
\[
g(t)=\langle \omega'(r_e+t(u-r_e))+e,u-r_e\rangle.
\]
It follows that $f$ is continuously differentiable on
$[0,1]$ with the derivative $g(t)$. Since the function attains its
minimum on $[0,1]$ at $t=0$, we have $g(0)\geq0$, which is exactly
(\ref{neweq12}).
\par
{\bf 2$^0$.}
At least the first statement of the following Lemma is well-known:
\begin{lemma}\label{lem1} For every $z\in Z^o$, the mapping $\xi\mapsto P_z(\xi)$
is a single-valued mapping of $\cE$ onto $Z^o$, and this mapping is
Lipschitz continuous, specifically,
\be
\|P_z(\zeta)-P_z(\eta)\|\leq \|\zeta-\eta\|_*\quad
\forall \zeta,\eta\in \cE. \ee{eq1}
Besides this, for all $u\in Z$, \be
\begin{array}{llll}
~\;\;&(a)& V(P_z(\zeta),u)&\leq V(z,u)+\langle
\zeta,u-P_z(\zeta)\rangle-V_z(z,P_z(\zeta))\\
&(b)&&\le V(z,u)+\langle \zeta,u-z\rangle+{\|\zeta\|_*^2\over
2}.
\end{array}
\ee{eq2}
\end{lemma}
\pr
Let $v\in P_z(\zeta)$, $w\in P_z(\eta)$. As
$V'_u(z,u)=\omega'(u)-\omega'(z)$, invoking (\ref{neweq12}), we have
$v,w\in Z^o$ and \be \langle
\omega'(v)-\omega'(z)+\zeta,v-u\rangle\leq 0\quad \forall u\in Z.
\ee{eq3} \be \langle \omega'(w)-\omega'(z)+\eta,w-u\rangle\leq
0\quad \forall u\in Z. \ee{eq3a} Setting $u=w$ in \rf{eq3} and
$u=v$ in \rf{eq3a}, we get
$$
\langle \omega'(v)-\omega'(z)+\zeta,v-w\rangle\leq 0,\,\, \langle
\omega'(w)-\omega'(z)+\eta,v-w\rangle\geq 0,
$$
whence $\langle
\omega'(w)-\omega'(v)+[\eta-\zeta],v-w\rangle\geq0$, or
$$\|\eta-\zeta\|_*\|v-w\|\geq \langle \eta-\zeta,v-w\rangle\geq
\langle \omega'(v)-\omega'(w),v-w\rangle\geq \|v-w\|^2,$$
and \rf{eq1} follows. This relation, as a byproduct, implies that
$P_z(\cdot)$ is single-valued.
\par
To prove \rf{eq2}, let $v=P_z(\zeta)$. We have
$$
\begin{array}{l} V(v,u)-V(z,u)\\
=[\omega(u)-\langle\omega'(v),u-v\rangle - \omega(v)]
-[\omega(u)-\langle \omega'(z),u-z\rangle - \omega(z)]\\
=\langle\omega'(v)-\omega'(z)+\zeta,v-u\rangle+\langle\zeta,u-v\rangle
-[\omega(v)-\langle\omega'(z),v-z\rangle-\omega(z)]\\
\leq\langle\zeta,u-v\rangle-V(z,v)\mbox{\ (due to \rf{eq3}),}\\
\end{array}
$$ as required in (\ref{eq:eq2}.$a$). The bound (\ref{eq:eq2}.$b$) is
obtained from (\ref{eq:eq2}.a) using the Young inequality:
\[
\langle \zeta, z-v\rangle\le {\|\zeta\|_*^2\over 2}+{1\over
2}\|z-v\|^2.
\]
Indeed, observe that by definition, $V(z,\cdot)$ is strongly
convex modulus 1 w.r.t. $\|\cdot\|$, and $V(z,v)\ge {1\over
2}\|z-v\|^2$, so that
\[
\langle \zeta, u-v\rangle-V(z,v)=\langle \zeta, u-z\rangle+\langle
\zeta, z-v\rangle-V(z,v)\le \langle \zeta,
u-z\rangle+{\|\zeta\|_*^2\over 2}.\eqno{\hbox{\epr}}
\]
\par
{\bf 3$^0$.}
We have the following simple corollary of Lemma \ref{lem1}:
\begin{corollary}\label{lem1-2}
Let $\xi_1,\xi_2,...$ be a sequence of elements of $\cE$. Define the
sequence $\{y_\tau\}_{\tau=0}^\infty$ in $Z^o$ as follows:
\[y_\tau=P_{y_{\tau-1}}(\xi_\tau),\;\;\;y_0\in Z^o.
\]
Then   $y_\tau$ is a measurable function of $y_0$ and
$\xi_1,...,\xi_\tau$ such that
\begin{equation}\label{aneq134}
(\forall u\in Z):\;\;\;\langle
-{\sum}_{\tau=1}^t \xi_\tau,u\rangle\le V(y_0,u)+{\sum}_{\tau=1}^t
\zeta_\tau,
\end{equation} with $|\zeta_\tau|\le
r\|\xi_\tau\|_*$ (here $r=\max_{u\in
Z}\|u\|$). Further,
\begin{equation}\label{aneq135}
{\sum}_{\tau=1}^t\zeta_\tau\le -{\sum}_{\tau=1}^t\langle\xi_\tau,y_{\tau-1}\rangle+
{1\over 2}{\sum}_{\tau=1}^t\|\xi_\tau\|^2_*.
\end{equation}
\end{corollary}
\pr
Using the bound (\ref{eq:eq2}.a) with $\zeta=\xi_\tau$ and
$z=y_{\tau-1}$, so that $y_\tau=P_{y_{\tau-1}}(\xi_\tau)$, we obtain for any
$u\in Z$:
\begin{equation}\label{aneq137} V(y_\tau,u)-V(y_{\tau-1},u)-\langle \xi_\tau,u\rangle\le
-\langle \xi_\tau,y_{\tau}\rangle-V(y_{\tau-1},y_\tau)\equiv \zeta_\tau; \end{equation}
summing up these inequalities over $\tau$ we get (\ref{aneq134}). Further, by definition of $P_z(\xi)$ we have
\[
\zeta_\tau=\max_{v\in Z}[-\langle \xi_\tau,v\rangle-V(y_{\tau-1},v)],
\]
so that $\zeta_\tau\leq r\|\xi_\tau\|_*$ due to $V\geq0$, and
\[
-r\|\xi_\tau\|_*\leq -\langle \xi_\tau,y_{\tau-1}\rangle =
[-\langle \xi_\tau,y_{\tau-1}\rangle-V(y_{\tau-1},y_{\tau-1})]\leq\zeta_\tau
\] due to $V(y_{\tau-1},y_{\tau-1})=0$. Thus, $|\zeta_\tau|\leq r\|\xi_\tau\|_*$, as claimed.
Further, by (\ref{eq:eq2}.b), where one should set $\zeta=\xi_\tau$, $z=y_{\tau-1}$, $u=y_\tau$, we have
\[
\zeta_\tau\le -\langle
\xi_\tau,y_{\tau-1}\rangle+{\|\xi_\tau\|_*^2\over 2}.
\]
Summing up these inequalities over $\tau$, we get (\ref{aneq135}).
\epr
\par{\bf 4$^0$.}
We also need the following result.
\begin{lemma}\label{lem:secret} Let   $z\in Z^o$, let
$\zeta$, $\eta$ be two points from $\cE$, and let
\[
w=P_z(\zeta),\qquad r_+=P_z(\eta)
\]
Then for all $u\in Z$ one has \be
\begin{array}{ll}
(a)&\|w-r_+\|\leq \|\zeta-\eta\|_*\\
(b) &V(r_+,u)-V(z,u)\le\langle \eta,u-w\rangle +\left[ \langle \eta,w-r_+\rangle-V(z,r_+)\right]\\
&\leq \langle
\eta,\,u-w\rangle+{1\over 2}\|\zeta-\eta\|_*^2-{1\over
2}\|w-z\|^2.
\end{array}
\ee{seq1}
\end{lemma}
\pr $(a)$: this is nothing but \rf{eq1}.
\par $(b)$: Using  (\ref{eq:eq2}.a) in Lemma \ref{lem1} we can write for $u=r_+$:
\[
V(w,r_+)\le V(z,r_+)+\langle \zeta,r_+-w\rangle-V(z,w).
\]
This results in \be V(z,r_+)\ge V(w,r_+)+V(z,w) +\langle
\zeta,w-r_+\rangle. \ee{vzz}
Now using (\ref{eq:eq2}.a)  with $\eta$
substituted for $\zeta$ we get
$$
\begin{array}{l}
V(r_+,u)\le V(z,u)+\langle \eta,u-r_+\rangle-V(z,r_+)\\
= V(z,u)+\langle \eta,u-w\rangle+\langle \eta,w-r_+\rangle-V(z,r_+)\\
\le V(z,u)+\langle \eta,u-w\rangle
+\langle \eta-\zeta,w-r_+\rangle
-V(z,w)-V(w,r_+)\mbox{\ [by \rf{vzz}]}\\
\le V(z,u)+\langle \eta,u-w\rangle +\langle
\eta-\zeta,w-r_+\rangle-{1\over 2}[\|w-z\|^2+\|w-r_+\|^2],
\\
\end{array}
$$
where the concluding inequalities are
due to the strong convexity of $\omega(\cdot)$. To conclude the bound $(b)$ of
\rf{seq1} it suffices to note that by the Young inequality,
\[
\langle \eta-\zeta,w-r_+\rangle\le {\|\eta-\zeta\|_*^2\over
2}+{1\over 2}\|w-r_+\|^2.\eqno{\hbox{\epr}}
\]
 \par{\bf 5$^0$.}
 We are able now to prove Theorem \ref{prop1prop1}.
By \rf{F} we have
\be
\|\wh{F}(w_\tau)-\wh{F}(r_{\tau-1})\|_*^2&\le&
(L\|r_{\tau-1}-w_\tau\|+M+\e_{r_{\tau-1}}+\e_{w_\tau})^2\nn &\le&
3L^2\|w_\tau-r_{\tau-1}\|^2+
3M^2+3(\e_{r_{\tau-1}}+\e_{w_\tau})^2.
\ee{Fb} Applying
Lemma \ref{lem:secret} with $z=r_{\tau-1}$,
$\zeta=\g_\tau\widehat{F}(r_{\tau-1})$,
$\eta=\g_\tau\widehat{F}(w_\tau)$ (so that $w=w_\tau$ and
$r_+=r_\tau$), we have for any $u\in Z$
$$
\begin{array}{l}
\langle \g_\tau \wh{F}(w_\tau),w_\tau-u\rangle+V(r_\tau,u)-V(r_{\tau-1},u)
\\
\le {\g_\tau^2\over 2}
\|\wh{F}(w_\tau)-\wh{F}(r_{\tau-1})\|^2-{1\over 2}\|w_\tau-r_{\tau-1}\|^2\\
\le {3\g_\tau^2\over 2}\left[L^2\|w_\tau-r_{\tau-1}\|^2+M^2+
(\e_{r_{\tau-1}}+\e_{w_\tau})^2\right]-
{1\over 2}\|w_\tau-r_{\tau-1}\|^2\hbox{\ [by \rf{Fb}]}\\
\le {3\g_\tau^2\over2}\left[M^2+(\e_{r_{\tau-1}}+\e_{w_\tau})^2\right] \hbox{\ [by \rf{gamma}]}\\
\end{array}
$$
When summing up from $\tau=1$ to $\tau=t$ we obtain
$$
\begin{array}{l}
{\sum}_{\tau=1}^{t} \langle \g_\tau
\wh{F}(w_\tau),w_\tau-u\rangle\\
\le
V(r_{0},u)-V(r_t,u)+{\sum}_{\tau=1}^t{3\g_\tau^2\over
2}\left[M^2+ (\e_{r_{\tau-1}}+\e_{w_\tau})^2\right]\\
\le
\Theta(r_0)+{\sum}_{\tau=1}^t{3\g_\tau^2\over 2}\left[M^2+
(\e_{r_{\tau-1}}+\e_{w_\tau})^2\right].\\
\end{array}
$$
Hence, for all $u\in Z$,
\begin{equation}\label{uu3001}
\begin{array}{l}
{\sum}_{\tau=1}^{t}
\langle \g_\tau F(w_\tau),w_\tau-u\rangle\\
\le \Theta(r_0)+{\sum}_{\tau=1}^t{3\g_\tau^2\over 2}
\left[M^2+(\e_{r_{\tau-1}}+\e_{w_\tau})^2\right]+{\sum}_{\tau=1}^t
\langle
\g_\tau\Delta_\tau,w_\tau-u\rangle\\
=\Theta(r_0)+{\sum}_{\tau=1}^t{3\g_\tau^2\over 2}
\left[M^2+(\e_{r_{\tau-1}}+\e_{w_\tau})^2\right]+{\sum}_{\tau=1}^t
\langle
\g_\tau\Delta_\tau,w_\tau-y_{\tau-1}\rangle\\
\quad+{\sum}_{\tau=1}^t\langle
\g_\tau\Delta_\tau,y_{\tau-1}-u\rangle\\
\end{array}
\end{equation} where $y_\tau$ are given by (\ref{recurrence}). Since the sequences
$\{y_\tau\}$, $\{\xi_\tau=\g_\tau\Delta_\tau\}$ satisfy the
premise of Corollary \ref{lem1-2}, we have
$$
\begin{array}{c}
(\forall u\in Z):\;\;{\sum}_{\tau=1}^t \langle
\g_\tau\Delta_\tau,y_{\tau-1}-u\rangle\le
V(r_0,u)+{\sum}_{\tau=1}^t{\g_\tau^2\over 2}\|\Delta_\tau\|_*^2\\
\le \Theta(r_0)+{\sum}_{\tau=1}^t{\g_\tau^2\over 2}\e^2_{w_\tau},\\
\end{array}
$$
 and thus (\ref{uu3001}) implies that for any $u\in Z$
\begin{equation}\label{uu300}
 {\sum}_{\tau=1}^{t} \langle \g_\tau
F(w_\tau),w_\tau-u\rangle\le\Gamma(t)
\end{equation}
with $\Gamma(t)$ defined in (\ref{eq:Gt0}). To complete the proof of \rf{zneq200} in the general case, note that since $F$ is monotone,
(\ref{uu300}) implies that for all $u\in Z$,
\[
{\sum}_{\tau=1}^t\g_\tau\langle F(u),w_\tau-u\rangle \leq
\Gamma(t),
\]
 whence
$$
\forall (u\in Z): \langle
F(u),\widehat{z}_t-u\rangle\leq\left[{\sum}_{\tau=1}^t\gamma_\tau\right]^{-1}\Gamma(t).
$$
When taking the supremum over $u\in Z$, we arrive at \rf{zneq200}.
\par
In the case of a Nash v.i., setting
$w_\tau=(w_{\tau,1},...,w_{\tau,m})$ and $u=(u_1,...,u_m)$ and
recalling the origin of $F$, due to the convexity of $\phi_i(z_i,z^i)$ in $z_i$, for all $u\in Z$ we get from (\ref{uu300}):
$$
\sum_{\tau=1}^{t} \g_\tau{\sum}_{i=1}^m[\phi_i(w_\tau)-\phi_i(u_i,(w_\tau)^i)]\le
\sum_{\tau=1}^{t} \g_\tau\sum_{i=1}^m\langle
F^i(w_\tau),(w_{\tau})_i-u_i\rangle
\leq \Gamma(t).
$$
 Setting
$\phi(z)={\sum}_{i=1}^m\phi_i(z)$, we get
\[
{\sum}_{\tau=1}^t\g_\tau\left[\phi(w_\tau)-{\sum}_{i=1}^m\phi_i(u_i,(w_\tau)^i)\right]\leq\Gamma(t).
\]
Recalling that $\phi(\cdot)$ is convex and $\phi_i(u_i,\cdot)$ are
concave, $i=1,...,m$, the latter inequality implies that
$$
\left[{\sum}_{\tau=1}^t\g_\tau\right]\left[\phi(\widehat{z}_t)-{\sum}_{i=1}^m\phi_i(u_i,(\widehat{z}_t)^i)\right]\leq\Gamma(t),
$$
or, which is the same,
$$
{\sum}_{i=1}^m\left[\phi_i(\widehat{z}_t)-{\sum}_{i=1}^m
\phi_i(u_i,(\widehat{z}_t)^i)\right]\leq\left[{\sum}_{\tau=1}^t\g_\tau\right]^{-1}
\Gamma(t).
$$
This relation holds true for all $u=(u_1,...,u_m)\in Z$; taking
maximum of both sides in $u$, we get
$$
\errNash(\widehat{z}_t)\leq
\left[{\sum}_{\tau=1}^t\g_\tau\right]^{-1} \Gamma(t).
\eqno{\hbox{\epr}}
$$

\subsection{Proof of Theorem \ref{themain}}
In what follows, we use the notation from Theorem
\ref{prop1prop1}.  By this theorem, in the case of constant
stepsizes $\gamma_\tau\equiv\gamma$ we have
\be
\errvi(\widehat{z}_t)\leq\left[t\gamma\right]^{-1}\Gamma(t),
\ee{bytheorem}
where
\be
\Gamma(t)&=&\Omega^2+{3\gamma^2\over
2}\sum_{\tau=1}^t\left[M^2+(\e_{r_{\tau-1}}+\e_{w_\tau})^2+{\e^2_{w_\tau}\over3}\right]+\gamma\sum_{\tau=1}^t \langle
\Delta_\tau,w_\tau-y_{\tau-1}\rangle\nn
&\leq& \Omega^2+{7\gamma^2\over
2}\sum_{\tau=1}^t\left[M^2+\epsilon_{r_{\tau-1}}^2
+\epsilon_{w_\tau}^2\right] +\gamma\sum_{\tau=1}^t
\langle \Delta_\tau,w_\tau-y_{\tau-1}\rangle.
\ee{bytheorem1}
For a Nash v.i., $\errvi$ in this relation can be replaced with
$\errNash$.
\par
Let us suppose that the random vectors $\zeta_i$ are defined on the probability space $(\Omega, \cF, \bP)$. We define two nested families of $\sigma$-fields
$\cF_i=\sigma(r_0, \zeta_1,\zeta_2,\ldots,\zeta_{2i-1})$ and $\cG_i=\sigma(r_0, \zeta_1,\zeta_2,\ldots,\zeta_{2i})$, $i=1,2,...$, so that $\cF_1\subset ...\,\cG_{i-1}\subset \cF_i\subset \cG_i\subset...$.
Then by description of the algorithm $r_{\tau-1}$ is $\cG_{\tau-1}$-measurable and $w_\tau$ is $\cF_{\tau}$-measurable. Therefore $\epsilon_{r_{\tau-1}}$ is $\cF_{\tau}$-measurable, and
$\epsilon_{w_\tau}$ and $\Delta_\tau$ are $\cG_{\tau}$-measurable. We conclude that under Assumption I we have
\begin{equation}\label{AssI}
\begin{array}{l}
\;\;\bE\left\{\epsilon_{r_{\tau-1}}^2|\cG_{\tau-1}\right\}\leq \sigma^2,\;
\bE\left\{\epsilon_{w_\tau}^2|\cF_{\tau}\right\}\leq \sigma^2,\;
\|\bE\left\{\Delta_\tau|\cF_{\tau}\right\}\|_*\leq \mu,\\
\end{array}
\end{equation}
and under Assumption II, in addition,
\begin{equation}\label{AssII}
\begin{array}{l}
\bE\left\{\exp\{\epsilon_{r_{\tau-1}}^2\sigma^{-2}\}\big|\cG_{\tau-1}\right\}\leq
\exp\{1\},\\
\bE\left\{\exp\{\epsilon_{w_\tau}^2\sigma^{-2}\}\big| \cF_{\tau}\right\}\leq \exp\{1\}.\\
\end{array}
\end{equation}
Now, let
\[
{\Gamma_0(t)}={7\gamma^2\over
2}{\sum}_{\tau=1}^t\left[M^2+\epsilon_{r_{\tau-1}}^2
+\epsilon_{w_\tau}^2\right].
\]
We conclude by  (\ref{AssI}) that
\begin{equation}\label{Gamma0}
\bE\left\{\Gamma_0(t)\right\}\leq {7\gamma^2t\over 2}[M^2+2\sigma^2].
\end{equation}
Further, $y_{\tau-1}$ clearly is
$\cF_{\tau-1}$-measurable, whence $w_\tau-y_{\tau-1}$ is $\cF_{\tau}$-measurable. Therefore
\be
\bE\left\{\langle
\Delta_\tau,w_\tau-y_{\tau-1}\rangle\big| \cF_{\tau}\right\}&=& \langle
\bE\left\{\Delta_\tau|\cF_{\tau}\right\},w_\tau-y_{\tau-1}\rangle\nn
&\leq&\mu\|w_\tau-y_{\tau-1}\| \leq
2\mu\Omega,
\ee{Gammarest}
where the concluding inequality follows from the fact that $Z$ is
contained in the $\|\cdot\|$-ball of radius
$\Omega$ centered at $z_c$, see (\ref{eqwesee}).
From \rf{Gammarest} it follows that
$$
\bE\left\{\gamma{\sum}_{\tau=1}^t \langle
\Delta_\tau,w_\tau-y_{\tau-1}\rangle\right\}\leq 2\mu\gamma
t\Omega.
$$
Combining the latter relation,
\rf{bytheorem}, \rf{bytheorem1} and (\ref{Gamma0}), we arrive at (\ref{K0}). (i) is proved.
\par
To prove (ii), observe, first, that setting
$$
J_t={\sum}_{\tau=1}^t\left[\sigma^{-2}\epsilon_{r_{\tau-1}}^2+
\sigma^{-2}\epsilon_{w_\tau}^2\right],
$$
we get
\begin{equation}\label{Gamma0astepsize}
\Gamma_0(t)={7\gamma^2M^2t\over 2}+{7\gamma^2\sigma^2\over 2}J_t.
\end{equation}
At the same time, setting $\cH_j=\sigma(r_0,\xi_1,...,\xi_j\}$, we can write
$$
J_t={\sum}_{j=1}^{2t} \xi_j,
$$
where $\xi_j\geq0$ is $\cH_j$-measurable, and
$$
\bE\left\{\exp\{\xi_j\}|\cH_{j-1}\right\}\leq\exp\{1\},
$$
see (\ref{AssII}). It follows that
\begin{equation}\label{under}
\begin{array}{l}
\bE\left\{\exp\{{\sum}_{j=1}^{k+1}\xi_j\}\right\}=
\bE\left\{\bE
\left\{\exp\{{\sum}_{j=1}^k\xi_j\}\exp\{\xi_{k+1}\}\right\}\big| \cH_k\right\}\\
=\bE\left\{\exp\{{\sum}_{j=1}^k\xi_j\}
\bE\left\{\exp\{\xi_{k+1}\}\big|\cH_k\right\}\right\}\leq
\exp\{1\}\bE\left\{\exp\{{\sum}_{j=1}^k\xi_j\}\right\}.\\
\end{array}
\end{equation}
Whence $\bE[\exp\{J\}]\leq\exp\{2t\}$, and applying
the Tchebychev inequality, we get
$$
\forall \Lambda>0:\;\Prob\left\{J>2t+\Lambda
t\right\}\leq\exp\{-\Lambda t\}.
$$
Along with (\ref{Gamma0astepsize}) it implies that
\begin{equation}\label{ForGamma0}
\forall \Lambda\geq0: \Prob\left\{\Gamma_0(t)>{7\gamma^2t\over
2}[M^2+2\sigma^2] + \Lambda{7\gamma^2\sigma^2t\over
2}\right\}\leq\exp\{-\Lambda t\}.
\end{equation}
Let now
$\xi_\tau=\langle\Delta_\tau,w_\tau-y_{\tau-1}\rangle$. Recall
that $w_\tau-y_{\tau-1}$ is $\cF_{\tau}$-measurable. Besides this, we have seen that
$\|w_\tau-y_{\tau-1}\|\leq D\equiv2\Omega$. Taking
into account (\ref{AssII}) and \rf{Gammarest}, we get
\begin{equation}\label{wegetastepsize}
\begin{array}{ll}
(a)&\bE\left\{\xi_\tau|\cF_{\tau}\right\}\leq\rho\equiv\mu D,\\
(b)&\bE\left\{\exp\{\xi_\tau^2R^{-2}\}|\cF_{\tau}\right\}
\leq\exp\{1\},\;\;\mbox{with}\;R=\sigma D.\\
\end{array}
\end{equation}
Observe that $\exp\{x\}\leq x+\exp\{9x^2/16\}$ for all $x$. Thus
 (\ref{wegetastepsize}.$b$) implies for $0\leq s\leq {4\over 3R}$
\be
\bE\left\{\exp\{s\xi_\tau\}|\cF_{\tau}\right\}&\leq&
\bE\{s\xi_\tau|\cF_{\tau}\}+\bE \left\{\exp\left\{{9s^2\xi^2_\tau\over 16}\right\}|\cF_{\tau}\right\}
\nn &\le &s\rho +
\exp\left\{{9s^2R^2\over 16}\right\}
\leq\exp\left\{s\rho+{9s^2R^2\over 16}\right\}.
\ee{concl1}
Further, we have $s\xi_\tau\leq {3\over 8}s^2R^2+{2\over
3}\xi_\tau^2R^{-2}$, hence for all $s\ge 0$,
$$
\bE\left\{\exp\{s\xi_\tau\}|\cF_{\tau}\right\}\leq\exp\{{3s^2R^2/
8}\}\bE\left\{\exp\left\{{2\xi_\tau^2\over 3R^2}\right\}|\cF_{\tau}\right\}\leq \exp\left\{{3s^2R^2\over 8}+{2\over 3}\right\}.
$$
When $s\geq{4\over 3R}$, the latter quantity is $\leq \exp\{3s^2R^2/4\}$,
which combines with \rf{concl1} to imply that for $s\ge 0$,
\begin{equation}\label{concl2}
\bE\left\{\exp\{s\xi_\tau\}|\cF_{\tau}\right\}\leq
\exp\{s\rho+3s^2R^2/4\}.
\end{equation}
Acting  as in (\ref{under}), we derive from (\ref{concl2}) that
$$
s\geq 0\Rightarrow
\bE\left\{\exp\{s{\sum}_{\tau=1}^t\xi_\tau\}\right\}\leq
\exp\{st\rho+3s^2tR^2/4\},
$$
and by the Tchebychev inequality, for all $\Lambda>0$,
$$
\Prob\left\{{\sum}_{\tau=1}^t\xi_\tau>t\rho+\Lambda
R\sqrt{t}\right\}\leq\inf_{s\geq0}\exp\{3s^2tR^2/4-s\Lambda
R\sqrt{t}\}=\exp\{-\Lambda^2/3\}.
$$
Finally, we arrive at
\begin{equation}\label{arriveat}\;~\;\;
\Prob\left\{\gamma{\sum}_{\tau=1}^t\langle\Delta_\tau,w_\tau-y_{\tau-1}\rangle>
2\gamma\left[\mu t+\Lambda
\sigma\sqrt{t}\right]\Omega\right\}\\
\leq\exp\{-\Lambda^2/3\}.\\
\end{equation}
for all $\Lambda>0$.
Combining \rf{bytheorem}, \rf{bytheorem1}, (\ref{ForGamma0}) and
(\ref{arriveat}), we get (\ref{K1}). \epr
\subsection{Proof of Lemma \ref{summarized}}
\paragraph{Proof of {\rm (i)}} We clearly have $Z^o=X^o\times Y^o$, and
$\omega(\cdot)$ is indeed continuously differentiable on this set.
 Let $z=(x,y)$ and $z'=(x',y')$, $z,z'\in Z$. Then
$$
\begin{array}{l}
\langle\omega'(z)-\omega'(z'),z-z'\rangle\\
={1\over\Omega_x^2}\langle
\omega_x^\prime(x)-\omega_x^\prime(x'),x-x'\rangle+
{1\over\Omega_y^2}\langle \omega_y^\prime(y)-\omega_y^\prime(y'),y-y'\rangle\\
\geq {1\over\Omega_x^2}\|x-x'\|_x^2+{1\over
\Omega_y^2}\|y-y'\|_y^2\geq \|[x'-x;y'-y]\|^2.\\
\end{array}
$$
Thus, $\omega(\cdot)$ is strongly convex on $Z$, modulus
$1$, w.r.t. the norm $\|\cdot\|$. Further, the minimizer of
$\omega(\cdot)$ on $Z$ clearly is $z_{\rm c}=(x_{\rm c},y_{\rm c})$, and it is immediately seen that $\max\limits_{z\in Z}V(z_{\rm c},z)=1$,
 whence
$\Omega=\sqrt{2}$.

\paragraph{Proof of {\rm (ii)}} {\bf 1$^0$.} Let $z=(x,y)$ and $z'=(x',y')$ with
$z,z'\in Z$. Note that by assumption {\bf G} in Section \ref{assumptions} we have
\begin{equation}\label{andthus}
\|y'\|_y\leq2\Omega_y. \end{equation} \par Further, we have from
(\ref{associated}) $F(z')-F(z)=\left[\Delta_x;\Delta_y\right]$, where
\bse
\Delta_x&=&{\sum}_{\ell=1}^m[\phi_\ell^\prime(x')-\phi_\ell^\prime(x)]^*[\bA_\ell y'+b_\ell]+
{\sum}_{\ell=1}^m[\phi_\ell^\prime(x)]^* \bA_\ell [y'-y],\\
\Delta_y&=&-{\sum}_{\ell=1}^m \bA_\ell^*[\phi_\ell(x)-\phi_\ell(x')]
+\Phi_*^\prime(y')-\Phi_*^\prime(y).
\ese
We have {\small
$$
\begin{array}{l}
 \|\Delta_x\|_{x,*}\\
  =\max\limits_{{h\in\cX\,\,\|h\|_x\leq1}}\langle h,
\sum\limits_{\ell=1}^m\left[[\phi_\ell^\prime(x')-\phi_\ell^\prime(x)]^*[\bA_\ell y'+b_\ell]+[\phi_\ell^\prime(x)]^*
\bA_\ell [y'-y]\right]
\rangle_\cX\\
\leq
\sum\limits_{\ell=1}^m\left[\max\limits_{{h\in\cX\,\atop\|h\|_x\leq1}}\langle
h,[\phi_\ell^\prime(x')-\phi_\ell^\prime(x)]^*[\bA_\ell y'+b_\ell]\rangle_\cX+
\max\limits_{{h\in \cX,\atop\|h\|_x\leq 1}}\langle
h,[\phi_\ell^\prime(x)]^* \bA_\ell [y'-y]
\rangle_\cX\right]\\
=\sum\limits_{\ell=1}^m\left[\max\limits_{{h\in\cX\,\atop\|h\|_x\leq1}}\langle
[\phi_\ell^\prime(x')-\phi_\ell^\prime(x)]h,\bA_\ell y'+b_\ell\rangle_\cX+
\max\limits_{h\in \cX,\|h\|_x\leq 1}\langle[\phi_\ell^\prime(x)]h,
\bA_\ell [y'-y]
\rangle_\cX\right]\\
\leq
\sum\limits_{\ell=1}^m\bigg[\max\limits_{{h\in\cX\,\|h\|_x\leq1}}\|[\phi_\ell^\prime(x')-\phi_\ell^\prime(x)]h\|_{(\ell)}\|\bA_\ell
y'+b_\ell\|_{(\ell,*)}\\
\multicolumn{1}{r}{+ \max\limits_{{h\in\cX\,\|h\|_x\leq1}}\|\phi_\ell^\prime(x)h\|_{(\ell)}\|\bA_\ell [y'-y]\|_{(\ell,*)}\bigg].}\\
\end{array}
$$
}
Then by (\ref{eq880}),{\small
$$
\begin{array}{l}
\|\Delta_x\|_{x,*}\\
\leq \sum\limits_{\ell=1}^m\bigg[[L_x\|x-x'\|_x+M_x][\|\bA_\ell
y'\|_{(\ell,*)}+\|b_\ell\|_{(\ell,*)}]
+\Omega_xL_x\|\bA_\ell[y-y']\|_{(\ell,*)}\bigg]\\
=[L_x\|x-x'\|_x+M_x]\sum\limits_{\ell=1}^m[\|\bA_\ell
y'\|_{(\ell,*)}+\|b_\ell\|_{(\ell,*)}]
+\Omega_xL_x\sum\limits_{\ell=1}^m\|\bA_\ell[y-y']\|_{(\ell,*)}
\\
\leq[L_x\|x-x'\|_x+M_x][\cA\|y'\|_y+\cB]+\Omega_xL_x\cA\|y-y'\|_y,\\
\end{array}
$$
}
by definition of $\cA$ and $\cB$. Next, due to (\ref{andthus}) we get  by definition of $\|\cdot\|${\small
\bse
\|\Delta_x\|_{x,*}
&\leq&
[L_x\|x-x'\|_x+M_x][2\cA\Omega_y+\cB]+\Omega_xL_x\cA\|y-y'\|_y\\
&\leq&
[\Omega_xL_x\|z-z'\|+M_x][2\cA\Omega_y+\cB]+\Omega_xL_x\cA\Omega_y\|z-z'\|,
\ese}
which implies
\begin{equation}\label{bounda}
 \|\Delta_x\|_{x,*}\leq
 \left[3\cA \Omega_xL_x\Omega_y+\cB\right]\|z-z'\|+[2\cA\Omega_y+\cB]M_x.
 \end{equation}
Further,{\small
\bse
\|\Delta_y\|_{y,*}&=&\max_{\eta\in\cY,\|\eta\|_y\leq1}\langle
\eta,-{\sum}_{\ell=1}^m \bA_\ell^*[\phi_\ell(x)-\phi_\ell(x')]
+\Phi_*^\prime(y')-\Phi_*^\prime(y)\rangle_\cY\\
&\leq& \max_{\eta\in\cY,\|\eta\|_y\leq1}{\sum}_{\ell=1}^m\langle
\eta,\bA_\ell^*[\phi_\ell(x)-\phi_\ell(x')]\rangle_{\cY}+\|\Phi_*^\prime(y')-\Phi_*^\prime(y)\|_{y,*}\\
&=&\max_{\eta\in\cY,\|\eta\|_y\leq1}{\sum}_{\ell=1}^m\langle
\bA_\ell \eta,\phi_\ell(x)-\phi_\ell(x')\rangle_{\cE_\ell}+\|\Phi_*^\prime(y')-\Phi_*^\prime(y)\|_{y,*}\\
&\leq&\max_{\eta\in\cY,\|\eta\|_y\leq1}{\sum}_{\ell=1}^m\|\bA_\ell\eta\|_{(\ell,*)}\|\phi_\ell(x)-\phi_\ell(x')\|_{(\ell)}
+\|\Phi_*^\prime(y')-\Phi_*^\prime(y)\|_{y,*}\\
&\leq&\max_{\eta\in\cY,\|\eta\|_y\leq1}{\sum}_{\ell=1}^m\|\bA_\ell\eta\|_{(\ell,*)}\Omega_xL_x\|x-x'\|_x
+[L_y\|y-y'\|_y+M_y],\ese}
by  (\ref{eq880}.$b$) and
(\ref{suchthat17}).
Therefore
$$
\|\Delta_y\|_{y,*}
\leq\cA \Omega_xL_x\|x-x'\|_x+L_y\|y-y'\|_y+M_y,
$$
whence
\begin{equation}\label{boundb}
\|\Delta_y\|_{y,*}\leq [\cA \Omega_x^2L_x+\Omega_yL_y]\|z-z'\|+M_y.
\end{equation}
From (\ref{bounda}), (\ref{boundb})  it follows that
$$
\begin{array}{l}
\|F(z)-F(z')\|_*\leq \Omega_x\|\Delta_x\|_{x,*}+\Omega_y\|\Delta_y\|_{y,*}\\
\leq
\left[4\cA L_x\Omega_x^2\Omega_y+\Omega_x\cB+\Omega_y^2L_y\right]\|z-z'\|+\left[2\cA\Omega_x\Omega_y+\Omega_x\cB\right]M_x+\Omega_yM_y.
\end{array}
$$
We have justified (\ref{eq900operator})
\par {\bf 2$^0$.} Let us verify  (\ref{eq900oracle}).
The first relation in (\ref{eq900oracle}) is readily given by
(\ref{eq900}.$a$,$c$). Let us fix $z=(x,y)\in Z$ and $i$, and let
\begin{equation}\label{leteq1}\begin{array}{rcl}
\Delta&=&F(z)-\Xi(z,\zeta_i)\\
&=&[\underbrace{{\sum}_{\ell=1}^m
[\phi_\ell^\prime(x)-\bG_\ell(x,\zeta_i)]^*\overbrace{[\bA_\ell
y+b_\ell]}^{\psi_\ell}}_{\Delta_x};
\underbrace{-{\sum}_{\ell=1}^m \bA_\ell^*[\phi_\ell(x)-f_\ell(x,\zeta_i)]}_{\Delta_y}.\\
\end{array}
\end{equation}
As we have seen, \begin{equation}\label{itfollows1}
{\sum}_{\ell=1}^m\|\psi_\ell\|_{(\ell,*)}\leq 2\cA\Omega_y+\cB.
\end{equation}
Besides this, for $u_\ell\in\cE_\ell$ we have
\begin{equation}\label{besidesthis}
\begin{array}{l}
\|{\sum}_{\ell=1}^m\bA_\ell^*u_\ell\|_{y,*}=
\max\limits_{{\eta\in\cY,\,\|\eta\|_y\leq1}}
\langle{\sum}_{\ell=1}^m\bA_\ell^*u_\ell,\eta\rangle_\cY=
\max\limits_{{\eta\in\cY,\,\|\eta\|_y\leq1}}
\langle{\sum}_{\ell=1}^mu_\ell,\bA_\ell\eta\rangle_\cY\\
\leq\max\limits_{{\eta\in\cY,\,\|\eta\|_y\leq1}}\left[{\sum}_{1\leq\ell\leq
m}\|u_\ell\|_{(\ell)}\|\bA_\ell\eta\|_{(\ell,*)}\right]\\
\leq\max\limits_{{\eta\in\cY,\,\|\eta\|_y\leq1}}\left[\max_{1\leq\ell\leq
m}\|u_\ell\|_{(\ell)}\right]{\sum}_{1\leq\ell\leq
m}\|\bA_\ell\eta\|_{(\ell,*)}
= \cA\max_{1\leq\ell\leq m}\|u_\ell\|_{(\ell)}.\\
\end{array}
\end{equation}
Hence, setting $u_\ell=\phi_\ell(x)-f_\ell(x,\zeta_i)$ we obtain
\begin{equation}\label{itfollows2}
\|\Delta_y\|_{y,*}=\|{\sum}_{\ell=1}^m\bA_\ell^*[\phi_\ell(x)-f_\ell(x,\zeta_i)]\|_{y,*}\leq
\cA\underbrace{\max_{1\leq\ell\leq
m}\|\phi_\ell(x)-f_\ell(x,\zeta_i)\|_{(\ell)}}_{\xi=\xi(\zeta_i)}.
\end{equation}
Further,
\bse
\|\Delta_x\|_{x,*}&=&\max\limits_{{h\in\cX,\,\|h\|_x\leq1}}\langle
h,{\sum}_{\ell=1}^m
[\phi_\ell^\prime(x)-\bG_\ell(x,\zeta_i)]^*\psi_\ell\rangle_\cX\\&=&
\max\limits_{{h\in\cX,\,\|h\|_x\leq1}}
{\sum}_{\ell=1}^m\langle[\phi_\ell^\prime(x)-\bG_\ell(x,\zeta_i)]h,\psi_\ell\rangle_\cX\\
&\leq&\max\limits_{{h\in\cX,\,\|h\|_x\leq1}}{\sum}_{\ell=1}^m\|
[\phi_\ell^\prime(x)-\bG_\ell(x,\zeta_i)]h\|_{(\ell)}\|\psi_\ell\|_{(\ell,*)}\\
&\leq&{\sum}_{\ell=1}^m\underbrace{\max\limits_{{h\in\cX,\,\|h\|_x\leq1}}\|[\phi_\ell^\prime(x)-\bG_\ell(x,\zeta_i)]h\|_{(\ell)}}_{\xi_\ell=\xi_\ell(\zeta_i)}\underbrace{\|\psi_\ell\|_{(\ell,*)}}_{\rho_\ell}\\
\ese
Invoking (\ref{itfollows1}), we conclude that
\begin{equation}\label{itfollows3}
\|\Delta_x\|_{x,*}\leq{\sum}_{\ell=1}^m\rho_\ell \xi_\ell,
\end{equation}
where all $\rho_\ell\geq0$, ${\sum}_\ell\rho_\ell\leq2\cA\Omega_y+\cB$ and
\[
\xi_\ell=\xi_\ell(\zeta_i)=\max\limits_{{h\in\cX,\,\|h\|_x\leq1}}\|[\phi_\ell^\prime(x)-\bG_\ell(x,\zeta_i)]h\|_{(\ell)}\\
\]
 Denoting by $p^2(\eta)$ the second moment of a scalar
 random variable $\eta$, observe that $p(\cdot)$ is a norm on the
space of square summable random variables representable as
deterministic functions of $\zeta_i$, and that
$$
p(\xi)\leq \Omega_xM_x,\,p(\xi_\ell)\leq M_x
$$
by  (\ref{eq900}.$b$,$d$).  Now by (\ref{itfollows2}),
(\ref{itfollows3}),
\bse
\left[\bE\left\{\|\Delta\|_*^2\right\}\right]^{{1\over 2}}&=&
\left[\bE\left\{\Omega_x^2\|\Delta_x\|_{x,*}^2
+\Omega_y^2\|\Delta_y\|_{y,*}^2\right\}\right]^{{1\over 2}}\\
&\leq&
p\left(\Omega_x\|\Delta_x\|_{x,*}+\Omega_y\|\Delta_y\|_{y,*}\right)
\leq p\left(\Omega_x{\sum}_{\ell=1}^m\rho_\ell\xi_\ell
+\Omega_y\cA\xi\right)\\
&\leq& \Omega_x{\sum}_\ell\rho_\ell \max_\ell p(\xi_\ell)+\Omega_y\cA
p(\xi)\\
&\leq& \Omega_x[2\cA\Omega_y+\cB]M_x+\Omega_y\cA\Omega_xM_x,
\ese
and the latter quantity is $\leq M$, see
(\ref{eq900operator}). We have established the second relation in
(\ref{eq900oracle}).
\par
{\bf 3$^0$.} It remains to prove that in the case of
(\ref{eq900str}), relation (\ref{eq900oraclestr}) takes place. To
this end, one can repeat word by word the reasoning from item
2$^0$ with the function $p_e(\eta)=\inf\left\{t>0:
\bE\left\{\exp\{\eta^2/t^2\}\right\}\leq\exp\{1\}\right\}$ in the
role of $p(\eta)$. Note that similarly to $p(\cdot)$, $p_e(\cdot)$
is a norm on the space of random variables $\eta$  which are
deterministic functions of $\zeta_i$ and are such that
$p_e(\eta)<\infty$. \epr
\subsection{Proof of Lemma \ref{summarized1}}
Item (i) can be verified exactly as in the case of Lemma \ref{summarized}; the facts expressed in (i) depend solely on
the construction from Section \ref{subsetup} preceding the latter Lemma, and are independent of what are the setups for $X,\cX$ and $Y,\cY$.\par
Let us verify item (ii). Note that we are in the situation
\begin{equation}\label{norms}
\begin{array}{l}
\|(x,\,y)\|=\sqrt{\|x\|_1^2/(2\ln(n))+|y|_1^2/(4\ln(p^{(1)}))},\\
\|(\xi,\,\eta)\|_*=\sqrt{2\ln(n)\|\xi\|_\infty^2+
4\ln(p^{(1)})|\eta|_\infty^2}.\\
\end{array}
\end{equation}
For $z=(x,y),z'=(x',y')\in Z$ we have
$$
F(z)-F(z')=\big[\underbrace{[\Tr((y-y')A_1);...;\Tr((y-y')A_n)]}_{\Delta_x};
\underbrace{-{\sum}_{j=1}^n(x_j-x_j^\prime)A_j}_{\Delta_y}\big].
$$
whence
\bse
\|\Delta_x\|_\infty&\leq& |y-y'|_1\max_{1\leq j\leq n}|A_j|_\infty\leq 2\sqrt{\ln(p^{(1)})}A_\infty\|z-z'\|,\\
|\Delta_y|_\infty&\leq& \|x-x'\|_\infty \max_{1\leq j\leq n}|A_j|_\infty\leq \sqrt{2\ln(n)}A_\infty\|z-z'\|,
\ese
and
\[ \|(\Delta_x,\,\Delta_y)\|_*\leq 2\sqrt{2\ln(n)\ln(p^{(1)})}A_\infty\|z-z'\|,
\]
as required in (\ref{neweq900operator}). Further, relation (\ref{neweq900oracle}.$a$) is clear from the construction of $\Xi_k$. To prove (\ref{neweq900oracle}.$b$), observe that when $(x,y)\in Z$, we have $\|\Xi^x(x,y,\eta)\|_\infty\leq\leq A_\infty$, $|\Xi^y(x,y,\eta)|_\infty\leq A_\infty$  (see (\ref{andset17}), (\ref{andset18})), whence
\be
\|\Xi^x(x,y,\eta)-F^x(x,y)\|_\infty\leq 2A_\infty,\,\,|\Xi^y(x,y,\eta)-F^y(x,y)|_\infty\leq 2A_\infty
\ee{wehaveseen}
due to $F(x,y)=\bE_\eta\left\{\Xi(x,y,\eta)\right\}$,
Applying \cite[Theorem 2.1(iii), Example 3.2, Lemma 1]{LargeDev}, we derive from \rf{wehaveseen} that for every $(x,y)\in Z$  and every $i=1,2,...$ it holds
\bse
&&\bE\left\{\exp\{\|\Xi_k^x(x,y,\zeta_i)-F^x(x,y)\|_\infty^2/N_{k,x}^2\}\right\}
\leq\exp\{1\},\\
&&N_{k,x}=2A_\infty\left(2\exp\{1/2\}\sqrt{\ln (n)}+3\right)k^{-1/2}
\ese and
\bse
&&\bE\left\{\exp\{\|\Xi_k^y(x,y,\zeta_i)-F^y(x,y)\|_\infty^2/N_{k,y}^2\}\right\}\leq\exp\{1\},\\
&&N_{k,y}=2A_\infty\left(2\exp\{1/2\}\sqrt{\ln (p^{(1)})}+3\right)k^{-1/2}.
\ese
Combining the latter bounds with  (\ref{norms}) we get  (\ref{neweq900oracle}.$b$). \epr


\begin{thebibliography}{ccc}
\bibitem{azuma}
Azuma, K. Weighted sums of certain dependent random variables.
{\em T\"okuku Math. J.}, {\bf 19} (1967), 357-367.
\bibitem{NERML}
Ben-Tal, A., Nemirovski, A. ``Non-Euclidean restricted memory
level method for large-scale convex optimization'' -- {\sl Math.
Progr.} {\bf 102} (2005), 407--456.
\bibitem{Khach}
Grigoriadis, M.D.,   Khachiyan, L.G. A sublinear-time randomized approximation algorithm
for matrix games. {\em Operations Research Letters} {\bf 18} (1995), 53--58.

\bibitem{Shapiroetal}
Juditsky, A. Lan, G., Nemirovski, A., Shapiro, A., Stochastic
Approximation
 Approach to Stochastic Programming,  {SIAM Journal on Opt.} {\bf 19} (2009) 1574-1609.
 \bibitem{LargeDev}
 Juditsky, A., Nemirovski, A. (2008), Large Deviations of Vector-valued Martingales in 2-Smooth
Normed Spaces\\
E-print: http://www.optimization-online.org/DB\_HTML/2008/04/1947.html
\bibitem{Korp}
Korpelevich, G. ``Extrapolation gradient methods and relation to
Modified Lagrangeans'' {\em Ekonomika i Matematicheskie Metody},
{\bf 19} (1983), 694-703 (in Russian; English translation in {\em
Matekon}).
%
\bibitem{NemYu}
 Nemirovski, A., Yudin, D., {\em Problem complexity and method efficiency in Optimization} { J. Wiley
\& Sons} (1983).
\bibitem{MP}
A. Nemirovski, ``Prox-method with rate of convergence $O(1/t)$ for
variational inequalities with Lipschitz continuous monotone
operators and smooth convex-concave saddle point problems'' --
{\sl SIAM J. Optim.} {\bf 15} (2004), 229-251.
\bibitem{Lunemmo}
Lu, Z., Nemirovski A., Monteiro, R. ``Large-Scale Semidefinite
Programming via Saddle Point Mirror-Prox Algorithm'', {\em Math.
Progr.}, {\bf 109} (2007), 211-237.
\bibitem{NOR}
Nemirovski, A., Onn, S., Rothblum, U. (2007), ``Accuracy certificates for computational problems with convex structure'' --
submitted to {\sl Mathematics of Operations Research}. E-print: http://www.optimization-online.org/DB\_HTML/2007/04/1634.html
\bibitem{nes1}
Nesterov, Yu. ``Smooth minimization of non-smooth functions'',
{\em Math. Progr.}, {\bf 103} (2005), 127-152.
\bibitem{nes2}
Nesterov, Yu. ``Excessive gap technique in nonsmooth convex
minimization'', {\em SIAM J. Optim.}, {\bf 16} (2005), 235-249.
\bibitem{nes3}
Nesterov, Yu. ``Dual extrapolation and its applications to solving
variational inequalities and related problems'', {\em Math.
Progr.} {\bf 109} (2007) 319-344.
\bibitem{ronit} Rubinfeld, R., Sublinear time algorithms. -- Marta Sanz-Sol\'{e},
Javier Soria, Juan Luis Varona, Joan Verdera, Eds. {\sl
International Congress of Mathematicians, Madrid 2006}, Vol. III,
 1o95--11110. European Mathematical Society Publishing House, 2006.
\end{thebibliography}
\end{document}